\documentclass[journal]{IEEEtran}
\usepackage[latin1,utf8]{inputenx}
\usepackage[dvips]{graphicx}
\usepackage{subfigure}
\usepackage[cmex10]{amsmath}
\usepackage{array}
\usepackage[colorlinks=true,citecolor=blue,linkcolor=blue,urlcolor=blue,plainpages=false]{hyperref}
\usepackage{booktabs, tabularx, multirow, array}
\usepackage{nomencl}
\usepackage{multirow}
\usepackage{comment} 
\usepackage{amsmath}
\usepackage{amsfonts}
\usepackage{amssymb}
\usepackage{graphicx}
\usepackage{amsthm}
\usepackage{eurosym}
\usepackage{algorithm}
\usepackage{tabularx}
\usepackage{cite}
\usepackage{array}
\usepackage{graphicx}
\usepackage[cmex10]{amsmath}
\usepackage{amsfonts,amssymb,amsxtra,balance}
\usepackage{layouts}

\usepackage{enumitem}
 \usepackage[table,xcdraw]{xcolor}

\usepackage{epstopdf}
\usepackage{MnSymbol,wasysym}
\usepackage{algorithm}
\usepackage{algpseudocode}
\usepackage{enumitem}

\usepackage{multirow}
\usepackage{dcolumn}
\usepackage{color}

\usepackage{setspace}
\usepackage[export]{adjustbox}
\usepackage{adjustbox}
\usepackage{soul}
\usepackage{booktabs}

\usepackage{tikz}
\usetikzlibrary{shapes,arrows}

\usepackage[textsize=scriptsize,shadow]{todonotes}

\DeclareGraphicsExtensions{.pdf,.png,.jpg}
\graphicspath{{images/}}

\def\w{\omega}
\newcommand{\inmat}[1]{\mbox{\rm {#1}}}
\newcommand{\E}{\mathbb{E}}

\def\X{\textbf{X}}

\def\S{\mathcal{S}}

\def\bxi{\boldsymbol{\xi}}
\def\bmu{\boldsymbol{\mu}}
\def\balpha{\boldsymbol{\alpha}}
\def\bpi{\boldsymbol{\pi}}

\def\h{\mathbf{h}}
\def\x{\mathbf{x}}
\def\y{\mathbf{y}}
\def\f{\mathbf{f}}
\def\q{\mathbf{q}}
\def\bphi{\mathbf{\phi}}
\def\blambda{\mathbf{\lambda}}

\def\X{\mathcal{X}}

\def\W{\mathbf{W}}
\def\T{\mathbf{T}}
\def\B{\mathbf{B}}

\def\e{\mathbf{e}}

\usepackage{marginnote}
\makenomenclature

\begin{document}
\title{Distributionally Robust Transmission Expansion Planning: a Multi-scale Uncertainty Approach}

\author{Alexandre~Velloso,
	David~Pozo,~\IEEEmembership{Senior Member,~IEEE},
	Alexandre~Street,~\IEEEmembership{Senior Member,~IEEE}
	
		\thanks{The work of Alexandre Velloso was supported by Finep/PIPG -- Financiadora de Estudos e Projetos / Programa de Incentivo \`{a} P\'{o}s-Gradua\c{c}\~{a}o. The work of D. Pozo was supported by Skoltech NGP Program (Skoltech-MIT joint project). Alexandre Street would like to acknowledge the financial support from CNPq -- Conselho Nacional de Desenvolvimento Cient\'{i}fico e Tecnol\'{o}gico -- and FAPERJ -- Funda\c{c}\~{a}o de Amparo \`{a} Pesquisa do Estado do Rio de Janeiro. 
	}	}

\maketitle

\begin{abstract}

We present a distributionally robust optimization (DRO) approach for the transmission expansion planning problem, considering both long- and short-term uncertainties on the system demand and non-dispatchable renewable generation. On the long-term level, as is customary in industry applications, we address the deep uncertainties arising from social and economic transformations, political and environmental issues, and technology disruptions by using long-term scenarios devised by experts. In this setting, many exogenous long-term scenarios containing partial information about the random parameters, namely, the average and the support set, can be considered. For each long-term scenario, a conditional ambiguity set models the incomplete knowledge about the probability distribution of the uncertain parameters in the short-term operation. Consequently, the mathematical problem is formulated as a DRO model with multiple conditional ambiguity sets. The resulting infinite-dimensional problem is recast as an exact, although very large, finite mixed-integer linear programming problem. To circumvent scalability issues, we propose a new enhanced-column-and-constraint-generation (ECCG) decomposition approach with an additional Dantzig--Wolfe procedure. In comparison to existing methods, ECCG leads to a better representation of the recourse function and, consequently, tighter bounds. Numerical experiments based on the benchmark IEEE 118-bus system are reported to corroborate the effectiveness of the method.

\end{abstract}


\IEEEpeerreviewmaketitle \addcontentsline{toc}{section}{\nomname}

\section*{Nomenclature}

{   
This section lists the main notation used in this work. Additional symbols are explained in the text or are interpretable in the context using the following general rules. The symbols with superscript ``$(j)$" denote variables, sets or results corresponding to the $j$-th iteration of the solution method. 
The symbols with superscript ``$k$" denote variables, parameters, and results associated with the $k$-th extreme point (or scenario) of a given support set. The symbols with subscript ``$\w$" refer to variables, parameters or sets related to the long-term scenario $\w$. The symbols with subscript ``$t$" refer to variables, functions or sets related to the $t$-th time period. 

\subsection{Sets and Indices}
\begin{description}[labelindent=8pt ,labelwidth=47pt, labelsep=3pt, leftmargin =58pt, style =standard, itemindent=0pt]%
	
	\item[$\Xi, \Xi_{\w}$]{Support set of the random vector and conditional support set of the random vector under the long-term scenario $\w$.}
	\vspace{0.005cm} 
	
	\item[$\Omega$]{Set of long-term scenarios}.
\vspace{0.005cm} 
	
	\item[$\mathcal{D}_{\w}$]Conditional ambiguity set under the long-term scenario $\w$. 
\vspace{0.005cm} 
		
	\item[$\mathcal{E}_\w$]Set of the extreme points of $\Xi_{\w}$.
\vspace{0.005cm} 

	\item[$K_\w$]Set of indices of the extreme points of $\Xi_{\w}$.	
	\vspace{0.005cm} 
	\item[$\mathcal{P}_\w$]Set of probability measures conditioned on the long-term scenario $\w$.
	\vspace{0.005cm} 
	\item[$\S, \S_\w$]Sample space, and subset of the sample space associated with the long-term scenario $\omega$.
	\vspace{0.005cm} 
	\item[$\mathbb{S},\mathbb{S}_\omega$] Appropriate sigma-algebras for $\S$ and $\S_\omega$, respectively.
\vspace{0.005cm} 
	
	\item[$\mathcal{T}$]Set of time periods $t$.
\vspace{0.005cm} 
	\item[$\mathcal{X}$]Set of feasible investment plans.

\end{description}

\subsection{Functions}
\begin{description}[labelindent=8pt ,labelwidth=47pt, labelsep=3pt, leftmargin =58pt, style =standard, itemindent=0pt]%
	
	\item[$\tilde{\bxi}$]Measurable function (or random vector) modeling the uncertainty in the net demand. 
	\vspace{0.005cm} 
	\item[$\tilde{\bxi}_t$]Subvector of $\tilde{\bxi}$ related to period $t$.
	\vspace{0.005cm} 
	\item[$\tilde{\bxi}(s)$]Realization of $\tilde{\bxi}$ for scenario $s$.
	\vspace{0.005cm} 
	\item[$g(\x,\bxi,\w)$]Minimum-cost dispatch function for investment $\x$, realization $\bxi$, and long-term scenario $\w$.
	
 	\vspace{0.005cm}

	\item[$H_{DR}(\x,\w)$] {Distributionally robust} recourse function 
	for investment $\x$, under the long-term scenario $\w$.
	\item[$\overline{H}_{DR}(\x,\w)$]Upper bound associated with $H_{DR}(\x,\w)$.
	\vspace{0.005cm} 
	\item[$\underline{H}_{DR}(\x,\w)$]Lower bound associated with $H_{DR}(\x,\w)$.
	
	
	
	
	
	
	
\end{description}

\subsection{Constants and Parameters}

\begin{description}[labelindent=8pt ,labelwidth=47pt, labelsep=3pt, leftmargin =58pt, style =standard, itemindent=0pt]%

\item[$\epsilon,\varepsilon$]Tolerances for the inner and main loops (in monetary and percentage units, respectively).
\vspace{0.005cm} 

\item[$\blambda^{(-)}_{\w},\blambda^{(+)}_{\w}$]Vectors of imbalance costs for the long-term scenario $\w$.
\vspace{0.005cm} 
\item[$\underline{\bmu}_{\w},\overline{\bmu}_{\w}$]Lower and upper bounds for the expected value of $\tilde{\bxi}$ for the long-term scenario $\w$.
\vspace{0.005cm} 
\item[$\rho_\w$]Probability or multi-objective weight of the long-term scenario $\w$.
\vspace{0.005cm} 
\item[$\bxi$]Generic point of $\Xi$.
\vspace{0.005cm} 
\item[$\bxi_{\omega}^k$]$k$-th extreme point of $\Xi_\w$.
\vspace{0.005cm} 
\item[$\bxi_{\omega}^*$]{Extreme point of $\Xi_\w$ associated with the maximum reduced cost of $\underline{H}_{DR}(\x,\w)$.}


	\vspace{0.005cm} 
	\item[$\mathbf{A}$]Line-bus incidence matrix.
	\vspace{0.005cm} 
	
	{\item[$ \mathbf{b}_\w$] Vector of right-hand-side parameters of the operative model}.
	\vspace{0.005cm} 
	\item[$\B_\w$]Spatial decoupling matrix for $\bxi$ under the long-term scenario $\w$.
	 \vspace{0.005cm} 
	\item[$\B_{t,\w}$] Submatrix of $\B_\w$ for period $t$.
	\vspace{0.005cm} 
	\item[$\mathbf{C}$] {Auxiliary matrix for disjunctive constraints.}
		 
	
	\item[$\mathbf{c}_{inv}$]Vector of investment costs.
	\vspace{0.005cm} 
	\item[$\mathbf{c}_\w$]{Vector of generation costs for the thermal generators under the long-term scenario $\w$.}
	\vspace{0.005cm} 
	\vspace{0.005cm} 
	\item[$\overline{c}_\w^*$]{Maximum reduced cost for the problem related to $\underline{H}_{DR}(\x,\w)$.}
	\vspace{0.005cm} 
	\item[$d$]Dimension of the uncertainty vector.
	\vspace{0.005cm} 
	\item[$\mathbf{e}$] Vector of ones with appropriate dimension.
	\vspace{0.005cm} 

	\item[$\overline{\mathbf{F}}$, $\overline{\mathbf{f}}$]{Matrix and vector for transmission constraints.}
	
	\vspace{0.005cm} 
	\item[$\mathbf{G}$]Thermal generator-bus incidence matrix.
	
	\vspace{0.005cm} 
	\item[$\h_\w$]{Vector of objective function coefficients of the compact operative model.}
	
	\vspace{0.005cm} 
	
	\item[$L$]Maximum number of inner loop iterations.
		\vspace{0.005cm} 
	\item[$LB^{(j)}$]Lower bound of the main loop at iteration $j$.	
		\vspace{0.005cm} 
	\item[$M$]Maximum number of scenarios added to the master problem per iteration of the algorithm.

	\vspace{0.005cm} 
	\item[$\overline{\q}_\w$]{Vector of generation limits for the thermal units under the long-term scenario $\w$.}
	\vspace{0.005cm} 
	\item[$\mathbf{r}^{dw}_\w, \mathbf{r}^{up}_\w$] {Vectors of ramp-down and ramp-up limits between two consecutive periods for the thermal units under the long-term scenario $\w$.}
	\vspace{0.005cm} 
	\item[$\mathbf{S}$]{Angle-to-flow matrix.} 
	\vspace{0.005cm} 

    \vspace{0.005cm} 
	\item[$\W_\w, \T_\w$]Recourse matrix and technology matrix for the long-term scenario $\w$.
	\vspace{0.005cm} 
	
	\item[$UB^{(j)}$]Upper bound of the main loop at iteration $j$.	
	\vspace{0.005cm} 
	
	\item[$z^*_{DR}$]{Optimal value for the distributionally robust transmission expansion planning problem.}

\end{description}

\subsection{Decision Variables or Vectors}
\begin{description}[labelindent=8pt ,labelwidth=47pt, labelsep=3pt, leftmargin =58pt, style =standard, itemindent=0pt]%
	
	\item[$\alpha_{0\w}$]Dual variable associated with the sum-one probability constraint.
		\vspace{0.005cm} 
	\item[$\underline{\balpha}_{\w}, \overline{\balpha}_{\w}$]{Dual vectors associated with the lower and upper limits, respectively, for the expected value of $\tilde{\bxi}$ under long-term scenario $\w$. }
		\vspace{0.005cm} 
	\item[$\mathbf{\theta}_{t}$]Vector of phase angles in {period $t$.}
		\vspace{0.005cm} 
	\item[$\bphi^{(-)}_{t},\bphi^{(+)}_{t}$]Power imbalance vectors in {period $t$.}
		\vspace{0.005cm} 
	\item[$\bpi_t$] Vector of dual variables of the compact version of the minimum-cost dispatch model.	
		\vspace{0.005cm} 
	
	\item[$\f_t$]Vector of power flows in {period $t$.}
		\vspace{0.005cm} 
	{\item[$p^k$]Probability of the extreme point $k$.}
		\vspace{0.005cm} 
	\item[$\q_t$]Vector of thermal generation in period $t$.
	 
	\item[$\x$]Vector of investment decisions.
	\vspace{0.005cm} 
	\item[$\y$]Vector of variables for the compact version of the minimum-cost operational {problem}.


\end{description}

%
%
%
%
%
%
%
%
%
%
%

}

\section{Introduction}
\label{sec:Int}

\IEEEPARstart{T}{he} transmission expansion planning (TEP) problem is generally related to strategic policies, as the outcome of a transmission plan {extends} far beyond providing a simple least-cost connection between the generation and loads. {For example, it may directly or indirectly shape the economic development for covered regions, or even facilitate policies for fostering innovation in various generation technologies. {As for electrical aspects}, the system reliability, operational flexibility, reserves deliverability, and long-run adaptability \cite{Hobbs2016} are key concepts that are significantly affected by the selected transmission capacity updates.} {On} the uncertainty side, planners have been dealing with several deep uncertainties arising from social and economic transformations, political and environmental issues, and technology disruptions, among others. {In this context, the definition of coherent future scenarios is a necessary step for defining the TEP \cite{borjeson2006scenario}. }

{Practical TEP applications generally rely on demand growth and renewable integration scenarios. { These scenarios, hereinafter referred to as long-term scenarios, are projections that are based on hypotheses about long-term drivers. The definition of the drivers that are relevant for the construction of the long-term scenarios frequently involves stakeholders and experts in various fields. In terms of modeling, the projections of long-term scenarios are generally provided by econometric models with many explanatory variables representing the long-term drivers. In this setting, a long-term scenario is a description of the future state of the system; that is, a conditional expectation on a specific set of beliefs and hypotheses (priors) that were used to parameterize the drivers. Notwithstanding, in practice, system planners want to consider multiple long-term scenarios \cite{Pancho2014,moreira2018}. Therefore, different hypotheses about the drivers are conceived, where these hypotheses draw distinct consistent narratives of the economic, political, technological, and environmental factors.  \cite{borjeson2006scenario}. A coherent long-term scenario for demand growth and renewable integration is then built on each of these hypotheses.}

\subsection{Motivation and literature review}

Despite the ability to consider many long-term profiles for the unknown parameters, the traditional deterministic ``what-if'' approach{, labeled in this work as D-TEP,} may not be sufficient for addressing all layers of uncertainty in modern power systems. For example, the large increase in renewable generation (RG) has introduced new levels of intermittency and unpredictability to electrical systems. These new aspects have motivated a change in numerous paradigms for short-term operation as well as for the planning of transmission networks {(see \cite{moreira2016reliable, li2018robust} and references therein)}.

{In this context, it may be reasonable to} consider both the long- and short-term effects of the uncertainty. For instance, a long-term scenario with high integration of wind-related sources would exhibit significantly more short-term variability than another where a reduced expansion of the wind power technology occurs.  Thus, the representation of multi-scale uncertainty in TEP models has been receiving attention in the recent literature across various approaches (stochastic, robust, and distributionally robust). 
We refer the interest reader to {\cite{moreira2018}, \cite{li2018robust},} \cite{PozoPscc2018}, \cite{Liu2018}, and \cite{Conejo2018} for further details. 	

The most popular {framework} for decision-making under uncertainty is the stochastic optimization (SO) {approach}, which considers scenarios drawn  from { a} predefined probability distribution function (PDF) \cite{Pancho2014,Zhan2019}. {Unfortunately}, the {PDF of short-term uncertainty is} considered difficult to estimate, specially in long-term studies, where the market and system structure may experience deep changes. {This fact has motivated the use of the adaptive robust optimization (ARO) framework to address TEP problems, which we label as the ARO-TEP approach}. In ARO-TEP problems, uncertainties are, in general, represented by polyhedral uncertainty sets relying on few assumptions about the uncertainty factors (see \cite{moreira2016reliable} and \cite{Ruiz2015}). Under this framework, the solution is that which performs the best in the worst-case scenario.

It is worth mentioning that while part of the recent literature on ARO has considered probability agnostic models (see \cite{moreira2018}, \cite{moreira2016reliable}, and \cite{Ruiz2015}), relevant efforts have been made to account for the information extracted from data to devise more realistic descriptions of the short-term uncertainty. We refer the interested reader to \cite{Velloso2019,liu2016stochastic}, and \cite{BoZeng2015} for applications in short-term operational models and to \cite{Conejo2018} for a hybrid-robust-and-stochastic approach applied to the TEP problem. Nevertheless, in long-term TEP applications, the use of scenario-based approaches relying on current data may be questionable, as the structure of the uncertainties in the target period may significantly differ from that found via data \cite{Pancho2014}. 
 
As an alternative to {SO and ARO approaches}, the distributionally robust optimization (DRO) framework, first introduced in \cite{scarf1957min}, has been developed within a broad mathematical and operations research context in the 2010s (see \cite{goh2010distributionally,delage2010distributionally,Wiesemann2014,Kuhn2018} and references therein). Unlike the SO approach that requires full PDF specification, the assumption for DRO-based models is that only partial information about the distributions of the uncertain parameters is available. The DRO framework is a robust approach, which is based on the worst-case expected cost. This contrasts with the more conservative ARO framework, which considers the single worst-case scenario. Therefore, the DRO {framework produces} ambiguity-averse models; that is, models that are robust against {predefined sets} of probability distributions, the so-called ambiguity set{s}. 

More recently, DRO approaches have been applied to model uncertainty in power system problems such as congestion line management \cite{qiu2015distributionally}, economic dispatch \cite{wei2016}, security-constrained optimal power flow \cite{xie2018distributionally}, risk-based optimal power flow with dynamic line rating \cite{wang2018risk}, investment decisions in wind farms \cite{alismail2018optimal}, and unit commitment \cite{Gourtani2016,xiong2017distributionally,zhao2017distributionally}. Regarding the application of DRO to TEP (DRO-TEP), to the best of the authors' knowledge only a few works have been published so far. {The network security was addressed in \cite{alvarado2018transmission} and a data-driven ambiguity set was applied in \cite{Bagheri2017} to account for the estimation uncertainty of empirical probability distributions. In both works, notwithstanding the relevant contributions to the DRO-TEP literature, the multi-scale nature of the uncertainty is not addressed. }

Most of the aforementioned works under robust approaches (ARO and DRO) are characterized by the traditional two-stage decision framework  involving  first-stage decision, uncertainty realization, and second-stage decision. The second-stage decision is modeled by a \emph{recourse function}, which in power system applications frequently represents the system's redispatch or corrective actions after uncertainty realization. Even though they provide flexibility in decision-making under uncertainty, two-stage robust models are generally hard to solve, requiring decomposition techniques.       

The column-and-constraint-generation (CCG) algorithm, developed in \cite{BoZeng2011}, has been largely applied to tackle two-stage models in power system applications under both ARO and DRO frameworks. The CCG method is a decomposition technique that relies on an iterative process that alternates between a master problem and an oracle subproblem. In many applications, each block of constraints and variables of the master problem is associated with one possible realization (scenario) of the {uncertain} parameter. In summary, the master problem is a relaxed version of the problem's extensive formulation where only a small subset of scenarios (and their respective constraints and variables) is represented. The oracle subproblem is a search procedure that, for a given first-stage trial solution determined by the master {problem}, finds {the} scenario that, according to some metric, produces the highest cost (or the largest violation for a specific set of constraints) for the recourse function. The constraints and variables associated with the new scenario are then represented in the master problem, which becomes more constrained. 

The CCG method is tailored to the ARO framework, since the first-stage solution and the single scenario determined by the oracle subproblem are sufficient to compute, at each iteration, the exact worst-case value for the recourse function. However, in the DRO framework, the recourse function is not completely determined by a single scenario, as it deals with worst-case expected values. As a consequence, recent applications of the CCG algorithm to solve problems under the DRO framework have relied on loose upper bounds for the recourse function. Moreover, many iterations of the CCG algorithm may be required until all the relevant scenarios are determined and introduced to the master problem.} Hence, to further develop the DRO approach
, new tailored uncertainty models and new solution methods are needed.

\subsection{Contributions}
In this work, we extend the ideas and developments of \cite{PozoPscc2018} and present a new DRO-TEP model. We assume that the true joint distribution of the uncertain parameters is difficult to estimate. Nevertheless, according to industry practices, we assume that partial information regarding the random parameters is available. In the designed framework, such partial information; that is, the characterization of the conditional expected values and support sets for the short-term net demand, is extracted from the long-term scenarios, which are projected by experts \cite{borjeson2006scenario,Pancho2014,moreira2018}. It is relevant to emphasize that, in this setting, the full specification of the conditional probability distribution for the short-term uncertainty is not required and that multiple exogenous long-term scenarios can be considered. 
	
{In terms of modeling,} in order to couple long- and short-term uncertainties, we propose an extension of the traditional {two-stage DRO framework to consider multiple} ambiguity sets. This {extension requires the introduction of the concept} of conditional ambiguity set, which is formalized in this work. Therefore, current industry practices, involving the {consideration of multiple} long-term scenarios, are accommodated within the proposed {multiple} conditional ambiguity set parametrization for the DRO-TEP model\footnote{{The proposed extension finds parallels in recent advances in robust optimization. The Stochastic Robust Optimization model in \cite{liu2016stochastic} and the Extended Robust Model in \cite{BoZeng2015} propose extensions to the robust optimization framework to consider multiple uncertainty sets}.}.

As for the solution approach, we propose a decomposition algorithm tailored for DRO problems where the ambiguity set is parameterized by first moment information. The new method is referred to as the enhanced column-and-constraint-generation (ECCG) algorithm. The proposed scheme differs from the CCG method by an inner loop, which is based on a Dantzig-Wolfe procedure (DWP) that determines more than one scenario at each main iteration of the algorithm. As a result, the proposed ECCG method provides a better representation of the recourse function and consequently a tighter \emph{Dantzig-Wolfe-like} upper bound, which is also mathematically formalized. Finally, as opposed to the CCG algorithm, more than one scenario is included to the master problem at each main iteration, resulting in higher lower bounds. In order to determine the scenarios to be added to the master problem, we propose a procedure that ranks the new scenarios based on their contribution to the recourse function value and selects only those of best rank. The combination of tighter lower and upper bounds results in fewer iterations and reduced computational burden as compared to the benchmark CCG algorithm.

In summary, the contributions of this work are twofold:

\begin{enumerate}
	\item A new multi-scale DRO-TEP that is based on the concept of multiple conditional ambiguity sets. This is a novelty in the literature of DRO and suitable to current industry practices in TEP. 

	\item A new decomposition method (ECCG algorithm) that provides better approximations for the distributionally robust recourse function than the CCG algorithm, resulting in tighter bounds.  

\end{enumerate}

\subsection{Paper Organization}
\label{subsec.Organization}

The rest of the paper is organized as follows. In Section \ref{sec.UncerModl}, the multi-scale uncertainty modeling is described and formalized. The proposed DRO-TEP model is formulated in Section \ref{sec.Model}. The solution methodology (ECCG algorithm) is presented in Section \ref{sec.Methodology}. The case study and numerical experiments are reported in Section \ref{sec.Results}. Finally, conclusions are addressed in Section \ref{sec: conclusions}.

\section{Uncertainty modeling}
\label{sec.UncerModl}



{The proposed approach recognizes the need for a multi-scale representation of the uncertainties affecting the {TEP} problem. The importance of integrating long- and short-term uncertainties in TEP problems has been discussed recently in the literature \cite{Conejo2018}. In Section \ref{quali}, an illustrative example based on a qualitative description of the proposed uncertainty representation is provided, and in Section \ref{probdescription} {the mathematical formalization of the approach is presented}.

\subsection{Qualitative description of the uncertain parameters}
\label{quali} 

In this work, the load and RG uncertainty are decomposed into: 
\begin{enumerate}
	\item A long-term component, which unfolds along {many} years and represents the uncertainty on the expected demand growth and long-term renewable integration.
	\item A short-term component, which is characterized by the conditional variability of the net demand (demand minus RG) on an hourly scale.
\end{enumerate}

In regards to the long-term component, in general TEP applications, governmental institutions or private consulting companies derive a few coherent long-term scenarios\footnote{Scenario analysis is largely used in several real power systems such as WECC, Chile, ERCOT, CAISO, UK National Grid (see \cite{borjeson2006scenario} and \cite{moreira2018}).} with different structural hypotheses for explanatory variables characterizing the economic, political, technological, and environmental factors. Based on such hypotheses, the expected demand growth and renewable integration levels are projected in subsequent studies  (see \cite{moreira2018} for further details).

\begin{figure}[t!]
	\centering
	\includegraphics[width=1.0\columnwidth]{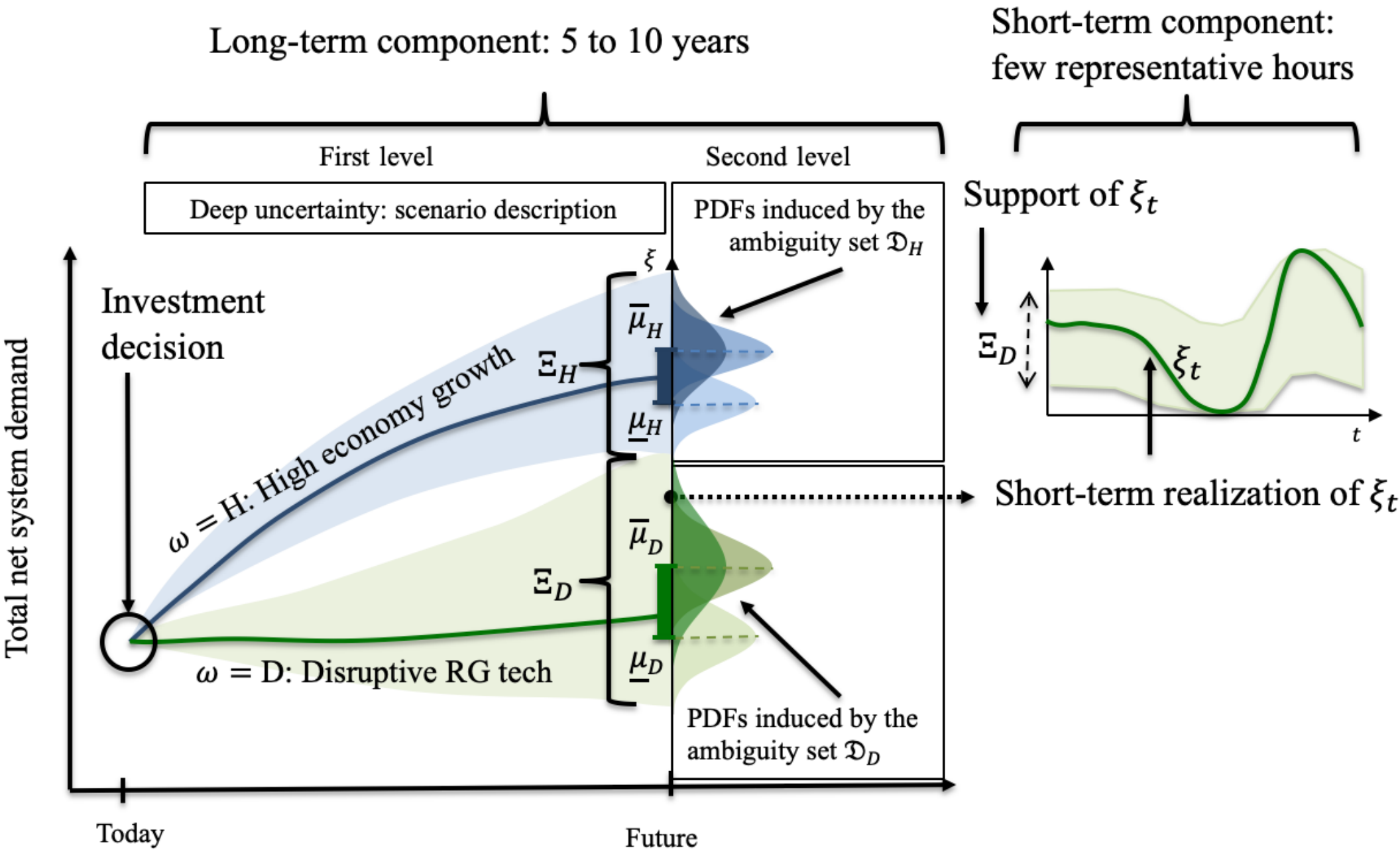}
	\caption{Uncertainty components of the net demand {(demand minus RG)}.}
	\label{fig: LongTermUncertainty}
\end{figure}

For example, experts may shape a plausible long-term scenario based on a future disruption in solar technology and the cost reduction of batteries and/or power electronic equipment. In this case, high penetration of distributed generation is expected, and therefore, a reduced expected system net demand (expected demand growth minus expected RG) is {projected}. However, the resulting short-term net load probability distribution may change significantly if a high economic development scenario with a lower distributed generation integration occurs. Thus, to allow the decision maker to use typical outcomes of long-term studies, in this work, a long-term scenario is characterized by two types of ranges: ranges for the expected value of the conditional probability distributions, and ranges for the possible outcomes of uncertainty factors (conditional support sets).

The short-term component of the uncertainty model is characterized by multiple conditional ambiguity sets (sets of conditional probability distributions), each of which connects to one corresponding long-term scenario. The information extracted from a long-term scenario to parameterize the related ambiguity set is the range for conditional expected value and range of the conditional support set. This setting allows planners to represent the lack of complete information regarding the underlying processes for the short-term net demand while considering partial information from long-term studies conducted by experts.

For comprehension purposes, a simple example, where the uncertain parameter is the total net demand, is illustrated in Fig. \ref{fig: LongTermUncertainty}. Two possible long-term scenarios are depicted, namely,  $\omega=H$ and $\omega=D$ representing ``\textbf{H}igh economic growth'' and ``\textbf{D}isruptive RG technology", respectively. The long-term scenarios determine the support sets ($\Xi_H$ and $\Xi_D$) and the ranges for the expected value of the net load  ($[\underline{\bmu}_{H}, \overline{\bmu}_{H}]$ and $[\underline{\bmu}_{D}, \overline{\bmu}_{D}]$). For each long-term scenario, the set of possible short-term distributions, induced by respective conditional ambiguity set, is abstractly represented in the illustration by different density profiles. In the proposed two-stage DRO framework, the investment decision (first-stage decision) is made under long- and short-term uncertainty components (represented by the circle on the intersection of the axes). The real-time operative dispatch decision (second-stage decision) is made under perfect information; that is, following the observation of both uncertainty components (signaled by the ``short-term realization" legend in the illustration). 
\subsection{Probabilistic description of the uncertain parameters}
\label{probdescription} 

{For the sake of simplicity, the only uncertainty source considered in this work is the net demand. It is represented by the vector of random factors $\tilde{\bxi}:\S\rightarrow\Xi$, with image set, or support set, $\Xi \subset \mathbb{R}^d$ defining all possible net demand extractions profiles. We consider a measurable space $(\S,\mathbb{S})$, where $\mathbb{S}$ is the appropriate sigma-algebra of the sample space $\S$, whose elements represent all possible states of nature. In this setting, $\tilde{\bxi}$ is a measurable function on $(\S,\mathbb{S})$ that maps the points from the sample space, $\S$, onto $\Xi$. The outcomes or realizations of $\tilde{\bxi}$ are represented by $\tilde{\bxi}(s)$ or, in short, $\bxi \in \Xi$, wherever convenient.

{In order to represent the concept of long-term scenarios} devised by experts, we consider $\Omega$ as the set of indices for the long-term scenarios, and we assume the existence of a partition, $\{\S_\w\}_{\w\in\Omega}$, of the sample space $\S$, i.e., $\cup_{\w\in\Omega}\S_\w = \S$ and $\S_\w \cap \S_{\w'}=\emptyset$. In this setting, conditional supports are the images of each part of the sample space, i.e., $\Xi_{\w}=\tilde{\bxi}(\S_\w)$, such that the unconditional support set is recovered by their union, i.e., $\cup_{\w\in\Omega} \Xi_{\w} = \Xi$. It is relevant to note that we do not assume that  $\Xi_\omega$ and $\Xi_{\omega'}$ are disjoint. 

The connection between the information provided by experts and the probabilistic framework is given by the { following inputs:}
\begin{enumerate} 
	\item[1)] Conditional expected value ranges, $\{[\underline{\bmu}_{\w}, \overline{\bmu}_{\w}]\}_{\w\in\Omega}$. 
	\item[2)] Conditional support sets, $\{\Xi_\w\}_{\w\in\Omega}$.
	
	\item[3)] Probabilities $\{\rho_\w=P(\S_\w)\}_{\w \in\Omega}$. 
	
	\item[3')] Alternatively, $\{\rho_\w\}_{\w \in\Omega}$ can also be defined as weights to generate Pareto-optimal solutions, thereby addressing long-term uncertainty under a multi-objective framework. 
\end{enumerate}

 }


Based on the conditional support set and expected value ranges for $\tilde{\bxi}$, we define the conditional ambiguity sets $\{\mathcal{D}_{\w}\}_{\w \in\Omega}$. Particularly, we assume that for each long-term scenario $\omega$, there exists a restricted (conditional) measurable space $(\S_\w \subset \S, \mathbb{S}_\w \subset \mathbb{S})$. The conditional ambiguity set $\mathcal{D}_{\w}$ is then defined as the set of all conditional probability measures (restricted to $\mathbb{S}_\w$) that induce a conditional expected value for $\tilde{\bxi}$ within the specified range, i.e.,
\begin{IEEEeqnarray}{l}
	\mathcal{D}_{\w} = \{P \in \mathcal{P}_{\w} : \underline{\bmu}_{\w} \leq   \E_P[\tilde{\bxi}\,|\,\S_{\w}] \leq \overline{\bmu}_{\w} \}. \label{eq. Def_D}
\end{IEEEeqnarray}

In \eqref{eq. Def_D}, $\mathcal{P}_\w$ represents the set of all probability measures in the measurable space $(\S_{\w}, \mathbb{S}_\w)$, and $\E_P[\tilde{\bxi}\,|\,\S_{\w}]$ represents the conditional expected value of $\tilde{\bxi}$ under a given probability measure $P$. For notation purposes, we denote by $\E_P[g(\tilde{\bxi})\,|\,\S_{\w}]:=\int_{\S_{\w}} g(\bxi) dP$ the conditional expectation of a given generic function $g(\tilde{\bxi})$ with regard to a measure $P\in \mathcal{D}_{\w}$ on $(\S_{\w}, \mathbb{S}_\w)$.

{For the sake of clarification, for the case depicted in Fig. \ref{fig: LongTermUncertainty}; that is, two long-term scenarios  $\Omega=\{\omega_1=H, \omega_2=D\}$, the inputs of the proposed uncertainty model would be: 
	1) the conditional expected value ranges, $[\underline{\bmu}_{H}, \overline{\bmu}_{H}]$ and $[\underline{\bmu}_{D}, \overline{\bmu}_{D}]$; 2) the conditional supports for the short-term component of the net load, $\Xi_H$ and $\Xi_D$; and 3) the probabilities (or weights) associated with each long-term scenario, $\rho_H$ and $\rho_D$.}


\section{Mathematical TEP model}
\label{sec.Model}

TEP {decisions (investments in candidate lines) are} represented by {the} binary vector $\x$. {The} set $\X$ defines the feasible investment plans. {On the operational side, the cost of operating the system under a given investment plan $\x$, net demand realization $\tilde{\bxi}(s)$, and long-term scenario $\w$ is represented by the minimum-cost dispatch function $g(\x,\tilde{\bxi}(s),\w)$.} Thus, in cases where the network planner {assumes} full knowledge regarding the probability measure {affecting the uncertain parameter $\tilde{\bxi}$}, i.e., under the assumption that $P(\cdot|\S_{\w})$ is the true conditional probability measure for each $\w\in\Omega$, the classical two-stage stochastic TEP approach is generally addressed by variants of the following optimization problem: 
\begin{IEEEeqnarray}{r'l'l}    \label{Pareto_Solution}
	\min_{\x \in \X}   &  \Big\{ \mathbf{c}_{inv}^{\top}\x + \sum_{\w\in\Omega} \rho_{\w} \E_P[g(\x,\tilde{\bxi},\w)|S_\w]  \label{stochastic}  \Big\}.    
\end{IEEEeqnarray}
Typically, the compact formulation for the operational cost function is given by  
\begin{IEEEeqnarray}{r'lll}
		g(\x,\bxi,\w)= \min _{\y\geq0} \Big\{\h^{\top}_\w\y  \Big| \,   \W_\w\y \geq \mathbf{b}_\w + \B_\w\bxi - \T_\w\x\Big\}, \quad  \label{compac}
\end{IEEEeqnarray}
whose right-hand-side is affected by an affine transformation of both investment and uncertainty vectors. Hence, $g(\x,\bxi,\w)$ is a convex function on $\x$ and $\bxi$. 

Note that the decision vector $\y$ stacks all the variables of the operational problem \eqref{compac}. The objective function of \eqref{compac} considers operational and imbalance costs through vector $\h_\w$, while the vectorial inequality comprises all operational constraints through matrices $\W_\w, \B_\w, \T_\w$, and vector $\mathbf{b}_\w$. 
	
In the next subsection (Section \ref{sec.TSO-TEP}), the compact formulation for the short-term operational model \eqref{compac} is detailed. The proposed DRO-TEP model is presented in Section \ref{sec.DRO-TEP} and a finite reformulation is derived in Section \ref{sec.finite}. In Section \ref{sec.particular}, we show how our model reduces to well-known benchmark models (ARO-TEP and D-TEP).   

\subsection{Short-term operation model description}
\label{sec.TSO-TEP}

The operational model is used to compute the dispatch cost under a given investment plan $\x\in\X$ and observed net load vector {$\bxi \in \Xi_{\omega}$, which can accommodate different components (sub-vectors), $\bxi_t$, per period.} Hence, the compact formulation presented in \eqref{compac} can be expanded as follows:
\begin{IEEEeqnarray}{ll}   
g(\x,\bxi,\w)	& =\hspace{-0.3cm} \min_{\substack{\q,\f,\theta\\ \bphi^{(-)}_{t},\bphi^{(+)}_{t} \ge \mathbf{0}}}   \sum_{t\in\mathcal{T}}( {\mathbf{c}^{\top}_\w\q_{t}} + \blambda^{(-)\top}_\w\bphi^{(-)}_{t} + \blambda^{(+)\top}_\w\bphi^{(+)}_{t}) \label{second_fobj} \quad \\
	 \quad &\inmat{s.t.:}    \nonumber \\
	  & \mathbf{G}\q_t +  \mathbf{A}\f_{t} = \B_{t,\w}\bxi_t - \bphi^{(-)}_{t}+\bphi^{(+)}_{t}, \; \forall t\in\mathcal{T} \label{second_balance}\vspace{+0.02cm}\\
	 &\big|\f_{t}  - \mathbf{S}\mathbf{\theta}_{t}\big| \leq \mathbf{C} (\e-\x),  \; \forall t\in\mathcal{T}    \label{second_line1}\vspace{-0.12cm}\\
	 &- \overline{\mathbf{F}}\x - \overline{\mathbf{f}} \leq \f_{t} \leq \overline{\mathbf{f}} + \overline{\mathbf{F}}\x, \; \forall t\in\mathcal{T} \label{second_line3} \\
	& -\mathbf{r}^{dw}_\w \leq \q_{t} - \q_{t-1}  \leq \mathbf{r}^{up}_\w, \; \forall t\in\mathcal{T} \label{eq.rampUP}\\
	& \mathbf{0} \leq \q_{t} \leq \overline{\q}_\w, \; \forall t\in\mathcal{T}. \label{eq.lastConst}
\end{IEEEeqnarray}

{ 
	
	As is customary in the literature of {TEP}, the standard optimal {dc} power flow model is used (see \cite{moreira2016reliable} and \cite{alvarado2018transmission}). The objective function \eqref{second_fobj} accounts {for the generation and imbalance costs} for all time periods.  
	{Constraint \eqref{second_balance} addresses the nodal power balance for all buses{. In this constraint, the matrix $\B_{t,\w}$ allocates, for each period $t$ and long-term scenario $\w$, the components of $\bxi_t$ (scenario representing uncertainty realization)} to the buses of the system.  
	{The relation \eqref{second_line1}} addresses the Kirchhoff's Voltage Law (KVL) through disjunctive  constraints, where the matrix $\mathbf{C}$ plays an auxiliary role. This matrix has as many rows as the total number of lines (existing and candidates) and as many columns as the number of candidate lines. The rows associated with the existing lines are composed of zeros, enforcing the KVL constraints. Each row representing a candidate line is composed of zeros and one nonzero component. Specifically, a large constant is assigned to the $i$-th column of the row representing the candidate line $i$. Thus, the KVL constraints are not enforced for the lines that are not selected. The inequalities in \eqref{second_line3} account for the transmission capacities of both existing and candidate lines. In these constraints,  $\overline{\mathbf{F}}$ is a diagonal matrix with the maximum capacity for the entries related to the candidate lines and zeros for the existing lines, while $\overline{\mathbf{f}}$ is a vector with the maximum capacity for existing lines and zeros for candidate lines. The expressions \eqref{eq.rampUP}--\eqref{eq.lastConst} represent the generation ramping limits and the maximum generation capacity respectively. Due to the imbalance terms in \eqref{second_balance}, problem \eqref{second_fobj}--\eqref{eq.lastConst} always admits a feasible solution.} }

\subsection{Distributionally robust TEP model}
\label{sec.DRO-TEP}

In order to address the lack of information regarding the true conditional probability measure $P(\cdot|\S_{\w})$, the recourse function is characterized by an ambiguity-averse preference functional, {based on the distributionally robust approach \cite{Wiesemann2014}}. Thus, the distributionally robust recourse function considers the worst-case conditional expected cost among all expectations induced by the probability measures in $\mathcal{D}_\w$. Mathematically, we have:
\begin{IEEEeqnarray}{l}
H_{DR}(\x,\w)=\sup_{P \in \mathcal{D}_{\w}}\mathbb{E}_P [g(\x,\tilde{\bxi},\w)|S_\w]. \label{inner}
\end{IEEEeqnarray}

Accordingly, the DRO-TEP problem is defined as an extension of \eqref{stochastic}:
\begin{IEEEeqnarray}{r'l'l}    \label{Pareto_DRO}
	z^*_{DR} = \min_{\x \in \X}   \Big\{ \mathbf{c}^{\top}_{inv}\x + \sum_{\w\in\Omega} \rho_{\w} H_{DR}(x,\w)    \Big\}.    
\end{IEEEeqnarray}

The $\sup$ problem in \eqref{inner} is the classical problem of moments \cite{landau1987moments}, which can be expressed as the following semi-infinite (infinitely many variables -- columns) linear program:

\begin{IEEEeqnarray}{rll"l}
  H_{DR}(\x,\w)=                &   \sup_{P \in \mathcal{P}_{\w}} &  \int_{\S_{\w}} g(\x,\tilde{\bxi},\w) dP                \label{def_H}   \\
	   \inmat{s.t.:} &  &\int_{\S_{\w}}  dP = 1 				&: \alpha_{0\w}  \quad \label{def_H1}  \\ 
	     	              & \underline{\bmu}_{\w} \leq 				& \int_{\S_{\w}} {\tilde{\bxi}} dP \leq \overline{\bmu}_{\w}   						&:\underline{\balpha}_{\w}, \,\overline{\balpha} _{\w}. \label{def_H2}
\end{IEEEeqnarray}
\subsection{Equivalent finite distributionally robust TEP model}
\label{sec.finite}
 By strong duality \cite{BoydBook}, problem \eqref{def_H}--\eqref{def_H2} admits the following semi-infinite (infinitely many constraints -- rows) dual formulation:
\begin{IEEEeqnarray}{l'll}
  H_{DR}(\x,\w)=   \min _{\alpha_0, \overline{\balpha}, \underline{\balpha}}    \left\{\alpha_{0\w} + \overline{\bmu}_{\w}^{\top} \overline{\balpha}_{\w} - \underline{\bmu}_{\w}^{\top}\underline{\balpha}_{\w}    \right\}   \label{H_infdual_eq0}  \\     
	    \inmat{s.t.:} 
	     \quad \alpha_{0\w}  +  \left(\overline{\balpha}_{\w} - \underline{\balpha}_{\w}\right)^{\top}\tilde{\bxi}(s) \ge  g(\x,\tilde{\bxi}(s),\w), \:  \forall s \in \S_{\w}.   \: \label{H_infdual_eq1}
\end{IEEEeqnarray}

The infinite set of linear constraints \eqref{H_infdual_eq1} can be recast as a single nonlinear worst-case constraint, \begin{IEEEeqnarray}{l'll}
	\alpha_{0\w} \ge \sup_{s \in \S_{\w}}\{ g(\x,\tilde{\bxi}(s),\w) - \left(\overline{\balpha}_{\w} - \underline{\balpha}_{\w}\right)^{\top}\tilde{\bxi}(s) \}. \label{nonlinConst}
\end{IEEEeqnarray}
This equivalence holds for the folllowing reason: Since \eqref{nonlinConst} holds for the supremum in $s$, it is also valid for all $s$ \cite{bental2009}. Note that the right-hand-side of \eqref{nonlinConst} is the supremum of a convex function on $\tilde{\bxi}(s)$ over a polyhedral set $\Xi_\w$. Therefore, the solution $s^*$ for \eqref{nonlinConst} is associated with a scenario whose image, $\tilde{\bxi}(s^*)$, belongs to the set of extreme points\footnote{This is a standard result from convex analysis -- the maximum of convex functions within a polyhedral set is achieved on one of the vertices \cite{BoydBook}.} of $\Xi_\w$, hereinafter referred to as $\mathcal{E}_\w=\{\bxi_{\omega}^k\}_{k\in K_\w}$. Consequently, if {constraint \eqref{H_infdual_eq1} is enforced for  $\mathcal{E}_\w$ only, the same} feasible region for $\alpha_0, \overline{\balpha}, \underline{\balpha}$ is found. Thus, problem \eqref{H_infdual_eq0}--\eqref{H_infdual_eq1} admits the following equivalent finite formulation:  
\begin{IEEEeqnarray}{l'll}
  H_{DR}(\x,\w)=   \min _{\alpha_0, \overline{\balpha}, \underline{\balpha}}    \left\{\alpha_{0\w} + \overline{\bmu}_{\w}^{\top} \overline{\balpha}_{\w} - \underline{\bmu}_{\w}^{\top}\underline{\balpha}_{\w}    \right\}   \label{H_dual_eq0}  \\     
	    \inmat{s.t.:} 
	     \quad \alpha_{0\w}  +  \left(\overline{\balpha}_{\w} - \underline{\balpha}_{\w}\right)^{\top}\bxi_\w^k \ge  g(\x,\bxi_\w^k,\w), \:  \forall k \in K_{\w},   \: \label{H_dual_eq1}
\end{IEEEeqnarray}

 \noindent which constitutes an equivalent finite formulation for the distributionally robust recourse function defined in \eqref{def_H}--\eqref{def_H2}. 
 
 { The equivalent finite primal formulation for the recourse function \eqref{def_H}--\eqref{def_H2} is achieved by computing the dual of problem  \eqref{H_dual_eq0}-\eqref{H_dual_eq1}}; that is,
\begin{IEEEeqnarray}{rrl'l}
	H_{DR}(\x,\w)=  & \max_{p^k} &  \sum_{k\in K_{\w}} g(\x,\bxi^k_{\omega},\w)\, p^k               \label{def_HMILP}   \\
	\inmat{s.t.:} &  &\sum_{k\in K_{\w}}  p^k = 1 & \;\;\qquad \label{def_HMILP1} \\ 
	& \underline{\bmu}_{\w} \leq 		& \sum_{k\in K_{\w}} {\bxi^k_{\omega}} \,p^k \leq \overline{\bmu}_{\w}.  \label{def_HMILP2} 
\end{IEEEeqnarray}
{As can be observed, the} candidates for worst-case distributions in $\mathcal{D}_{\w}$ are discrete.}

Aiming to achieve an equivalent finite formulation for the DRO-TEP model \eqref{Pareto_DRO}, it should be noted that the right-hand-side of \eqref{H_dual_eq1} is the optimal value of an instance of problem \eqref{second_fobj}--\eqref{eq.lastConst} for a given $\x$ and a scenario $\bxi_\w^k$ in $\mathcal{E}_\w$. For notational conciseness, we replace $g(\x,\bxi_\w^k,\w)$ in \eqref{H_dual_eq1} with the objective function of \eqref{compac}, namely $\h_\w^\top\y_\w^k$, evaluated on a generic \emph{feasible} operative point, $\y_\w^k$. In order to obtain an equivalent model, the operative variables $\y_\w^k$ and feasibility constraints of \eqref{compac} for scenario $\bxi_{\omega}^k$ must be included in the new formulation. This substitution results in the following equivalent extended dual formulation for the recourse function:
\begin{IEEEeqnarray}{l'll}
	H_{DR}(\x,\w)=   \min _{\y_\w^k, \alpha_0, \overline{\balpha}, \underline{\balpha}}    \left\{\alpha_{0\w} + \overline{\bmu}_{\w}^{\top} \overline{\balpha}_{\w} - \underline{\bmu}_{\w}^{\top}\underline{\balpha}_{\w}    \right\}   \label{H_dualEXT_eq0}  \\     
	\inmat{s.t.:} 
	\quad \alpha_{0\w}  +  \left(\overline{\balpha}_{\w} - \underline{\balpha}_{\w}\right)^{\top}\bxi_\w^k \ge  \h^\top_\w\y_\w^k, \:  \forall k \in K_{\w}   \: \label{H_dualEXT_eq1}\\
	\quad \quad \;\, \W_\w \y_\w^k\; \geq \, \mathbf{b}_\w + \B_\w \bxi_\w^k - \T_\w \x,   \forall k \in K_{\w}  \label{H_dual_operEXT_eq2}.
\end{IEEEeqnarray}

The resulting set of constraints in \eqref{H_dualEXT_eq1}--\eqref{H_dual_operEXT_eq2} is equivalent to \eqref{H_dual_eq1} (in the sense of producing the same optimal value for problem \eqref{H_dual_eq0}--\eqref{H_dual_eq1}), because 1) by optimality, $\h^\top\y_\w^k$ is bounded from below by $g(\x,\bxi_\w^k,\w)$; and 2) every point that is feasible for \eqref{H_dual_eq1} is also feasible for the new set of constraints. 

Thus, {the DRO-TEP problem \eqref{Pareto_DRO} can be recast as the following equivalent finite mixed-integer scenario-based linear program}:
\begin{IEEEeqnarray}{rll}
	z^*_{DR} & =  \min _{\substack{\x, \y_\w^k, \alpha_{0\w},\\ \overline{\balpha}_\w, \underline{\balpha}_\w }}  \mathbf{c}^{\top}_{inv}\x + \sum_{\w\in\Omega} \rho_\w \left[\alpha_{0\w} + \overline{\bmu}_{\w}^{\top} \overline{\balpha}_{\w} - \underline{\bmu}_{\w}^{\top}\underline{\balpha}_{\w} \right]  \label{H_dual_oper_eq0}   \quad \\     
		\inmat{s.t.:}  &  \quad \x \in \X \\
		& \h^\top_\w\y_\w^k \leq \alpha_{0\w}  +  \left(\overline{\balpha}_{\w} - \underline{\balpha}_{\w}\right)^{\top}\bxi_\w^k, \forall \w\in\Omega, k \in K_{\w}   \: \label{H_dual_oper_eq1}\\
		&  \W_\w \y_\w^k\; \geq \, \mathbf{b}_\w + \B_\w \bxi_\w^k - \T_\w \x,   \forall \w\in\Omega, k \in K_{\w}  \label{H_dual_oper_eq2}.
\end{IEEEeqnarray}
\subsection{Particular cases of the distributionally robust TEP model}
\label{sec.particular}

The DRO-TEP \eqref{H_dual_oper_eq0}--\eqref{H_dual_oper_eq2} can be particularized to the ARO-TEP approach. 
 The ARO-TEP formulation 
	can achieved by disregarding the moment information; that is, dropping the terms involving $\underline{\balpha}_{\w}$ and $\overline{\balpha}_{\w}$ of expressions \eqref{H_dual_oper_eq0} and \eqref{H_dual_oper_eq1}. Under such an assumption, for each  $\omega \in \Omega$, the Cartesian product of (each dimension of) the range  for the expected value $[\underline{\bmu}_{\w}, \overline{\bmu}_{\w}]$ is equivalent to the whole support set $\Xi_\w$. Thus, except for degenerate cases, the worst-case probability measure { for the ARO-TEP model} is such that all probability mass is assigned to the worst-case short-term scenario.

In the same manner, the DRO-TEP and ARO-TEP models can also be reduced to the D-TEP model. To do so, one just needs to collapse both the support set and the moment range to a singleton (scenario). The aforementioned particularizations result in two benchmarks models, D-TEP and ARO-TEP, one less and the other more conservative than the proposed DRO-TEP model, respectively. These benchmarks are used for comparison purposes in our computational experiments and do not represent all possible variants within their specific classes. 

\section{Solution methodology}
\label{sec.Methodology}

Problem \eqref{H_dual_oper_eq0}--\eqref{H_dual_oper_eq2} can be solved directly, in what we label as the \emph{full-vertex approach} (FVA) \cite{PozoPscc2018}. However, this approach is not scalable since the problem size increases exponentially with the dimension ($d$) of the uncertain parameter. Therefore, decomposition techniques are required.

 In Section \ref{subCCG}, we discuss the benchmark CCG algorithm and the issues related to its application to DRO-based problems. An overview of the proposed ECCG algorithm is presented in Section \ref{subOverviewECCG} and the detailed ECCG algorithm is reported in Section \ref{subDetailECCG}. The ECCG algorithm is general enough to be applied to similar problems (e.g. \cite{zhao2017distributionally}), however, for didactic purposes, the discussions in this section are focused on the proposed DRO-TEP model.


\subsection{Overview of the traditonal CCG algorithm}
\label{subCCG}

The CCG algorithm alternates between a master problem and an oracle subproblem. The master problem is a relaxed version of \eqref{H_dual_oper_eq0}--\eqref{H_dual_oper_eq2} considering a small subset of constraints \eqref{H_dual_oper_eq1}--\eqref{H_dual_oper_eq2}. We use $K_\w^{(j)}\subset K_\w$ to represent the subset of constraints at iteration $j$, which are associated with the subset of scenarios $\mathcal{E}_\w^{(j)}\subset \mathcal{E}_\w$ (vertices of $\Xi_\w$). At each iteration $j$, the master problem determines a lower bound ($LB^{(j)}$) for $z^*_{DR}$ and a trial TEP solution $\x^{(j)}$. The master problem also determines, for each $\omega \in \Omega$, the variables associated with the recourse function ($\alpha_{0\w}^{(j)},\overline{\balpha}_{\w}^{(j)},\underline{\balpha}_{\w}^{(j)}$) and the operational response ($\y_\w^{(j)}$). In turn, for each $\omega \in \Omega$, an oracle subproblem identifies the new scenario $\bxi_\w^*$ (not yet included to $\mathcal{E}_\w^{(j)}$) that most violates constraint \eqref{H_dual_oper_eq1}, which is associated with constraint \eqref{H_dual_eq1} in the dual recourse function problem. Equivalently, $\bxi_\w^*$ is the scenario that is related to the new primal variable (new column $p$), which features the maximum reduced cost, $\overline{c}_\w^*$, for the constrained primal recourse function problem; that is, \eqref{def_HMILP}--\eqref{def_HMILP2} restricted to the scenarios in $\mathcal{E}_\w^{(j)}$.

The violated constraints and variables associated with $\bxi_\w^*$ (for all $\omega \in \Omega$) are added to the master problem at the next iteration. The algorithm terminates when no violated constraints are found by the oracle subproblem or when the upper bound, $UB^{(j)}$, is arbitrarily close to the $LB^{(j)}$.

However, as opposed to the ARO framework, where the value of the worst-case recourse function is given by $g(\x,\bxi_\w^*,\w)$, the definition of $H_{DR}(\x^{(j)},\w)$ requires more than a single scenario. This implies that approximations that generally involve the reduced cost $\overline{c}_\w^*$ are necessary to compute an upper bound. Unfortunately, when the scenarios in $\mathcal{E}_\w^{(j)}$ are insufficient for a reasonable approximation of $H_{DR}(\x,\w)$, the values of $\overline{c}_\w^*$ may be very large, therefore generating loose upper bounds. Moreover, many iterations of the method may be required until all relevant scenarios are determined and included in the master problem.

\subsection{Overview of the ECCG algorithm}
\label{subOverviewECCG}

	\begin{figure}[t!]
		\centering
		\includegraphics[width=0.91\columnwidth]{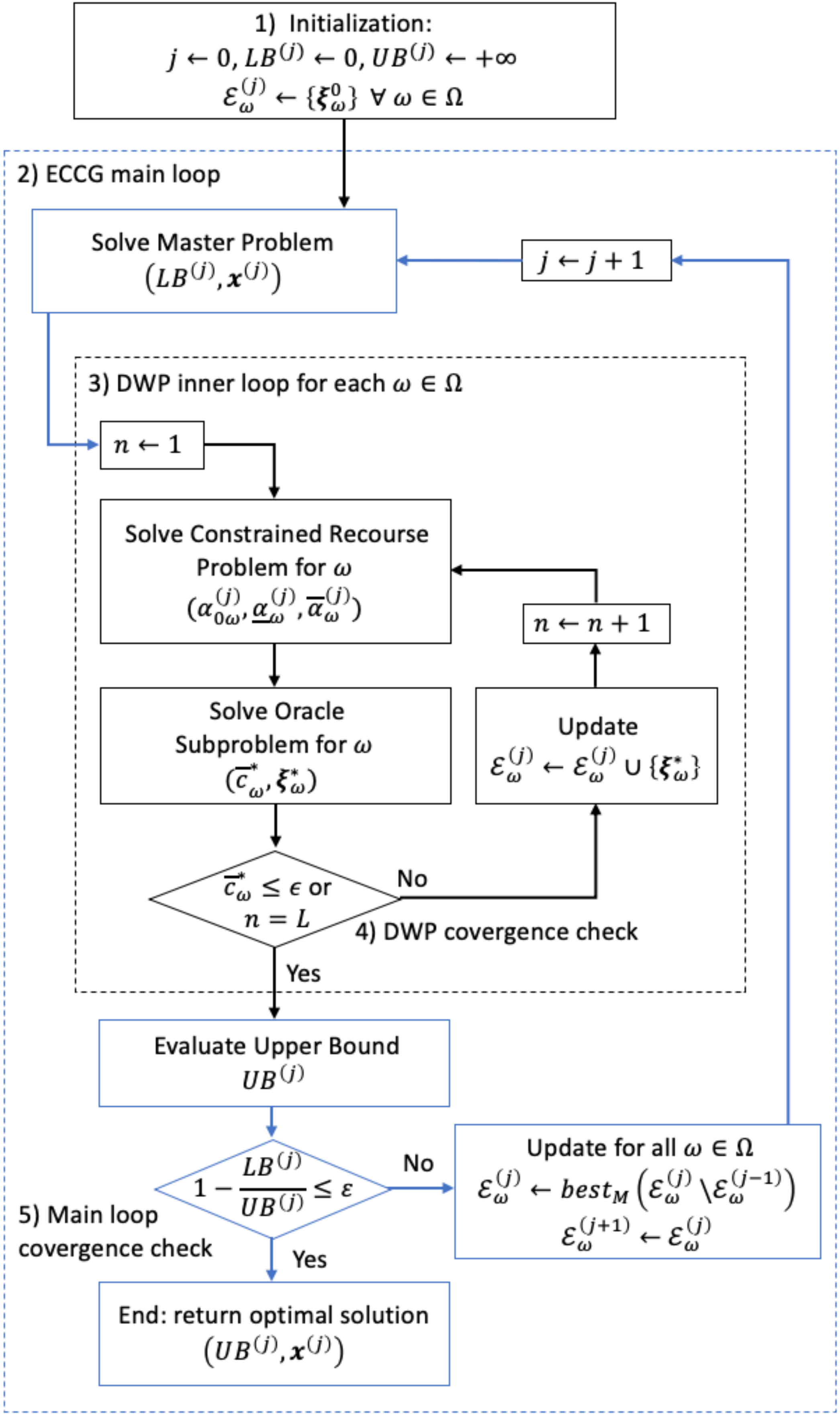}
		\caption{Flowchart of the proposed {ECCG algorithm}.}
		\label{fig: ECCG}
	\end{figure}
	
The determination of a tight upper bound $UB^{(j)}$ involves the computation of the true value of the recourse function for all $\w\in\Omega$; that is, $\{H_{DR}(\x^{(j)},\w)\}_{\w\in\Omega}$. Unfortunately, computing $H_{DR}(\x^{(j)},\w)$ might be intractable for high dimensional instances, since it requires the solution of problem \eqref{def_HMILP}--\eqref{def_HMILP2} or its dual, \eqref{H_dual_eq0}--\eqref{H_dual_eq1}, considering the entire set of exponentially many extreme points $\mathcal{E}_\w$. Thus, in order to deal with scalability issues while controlling the approximation of the recourse function, we have extended the CCG algorithm by introducing an inner loop. The inner loop is the DWP applied to problem \eqref{def_HMILP}--\eqref{def_HMILP2}.

At each main iteration $j$, the DWP performs up to $n=L$ inner iterations, where a constrained version of the recourse function problem \eqref{def_HMILP}--\eqref{def_HMILP2} (considering only scenarios in $\mathcal{E}_\w^{(j)}$) is alternated with the oracle subproblem. The constrained version of the recourse function problem updates ($\alpha_{0\w}^{(j)},\overline{\balpha}_{\w}^{(j)},\underline{\balpha}_{\w}^{(j)}$) and determines a lower (suboptimal) approximation for the true recourse function. The oracle subproblem finds a new short-term scenario $\bxi_\w^*$ (featuring the maximum reduced cost, $\overline{c}_\w^*$) which is then introduced to the constrained recourse problem \eqref{def_HMILP}--\eqref{def_HMILP2}. The DWP converges when the maximum reduced cost is below $\epsilon$ or the maximum number of (inner) iterations $L$ is reached. As is shown in Section \ref{teo1}, the upper bound for the ECCG algorithm is controllable and $\epsilon$-tight for large values of $L$.

The final step of the ECCG algorithm is the selection of the scenarios that are added to the master problem. Because the master problem is the bottleneck of the DRO-TEP problem, it is not advised to include all scenarios determined by the DWP at each main iteration. Thus, the following selection procedure is adopted. The new scenarios are ranked according to their contribution to the objective value of the constrained recourse problem in the last inner iteration. Then, the $M \le L$ scenarios with highest ranks are added to the master problem. For comprehension purposes, the proposed algorithm is summarized in Fig.\ref{fig: ECCG}. 


\subsection{Detailed description of the ECCG algorithm}
\label{subDetailECCG}
{The detailed algorithm is as follows:}

\vspace{0.2 cm}
\noindent \textbf{1) Initialization:} The set of scenarios $\mathcal{E}_\w^{(j)}$ is initialized at $j\leftarrow0$ with a dummy scenario $\bxi_{\omega}^{0}\leftarrow(\overline{\bmu}_{\w}+\underline{\bmu}_{\w})/2$ to ensure the existence of a probability measure capable of producing an expected value within $[\underline{\bmu}_{\w},\overline{\bmu}_{\w}]$; that is,  $\mathcal{E}_\w^{(0)}\leftarrow\{\bxi_{\omega}^{0}\}$. 

\vspace{0.2 cm}
\noindent \textbf{2) ECCG main loop (assessing $LB^{(j)}$ and $\x^{(j)}$):} For each iteration $j$ of the main loop, the master problem solution provides a trial solution $\x^{(j)}$ and its corresponding lower bound $LB^{(j)}$. The master problem is expressed as follows: 
\begin{IEEEeqnarray}{l'll}
	\min_{\substack{\x, \y_\w^k, \alpha_{0\w},\\ \overline{\balpha}_\w, \underline{\balpha}_\w }}  \mathbf{c}^{\top}_{inv}\x + \sum_{\w\in\Omega} \rho_\w \left[\alpha_{0\w} + \overline{\bmu}_{\w}^{\top} \overline{\balpha}_{\w} - \underline{\bmu}_{\w}^{\top}\underline{\balpha}_{\w} \right]     \label{Master1}  \\     
	\inmat{s.t.:}  
	\; \x \in \X \\
	\;\,\;\quad  \alpha_{0\w}  +  \left(\overline{\balpha}_{\w} - \underline{\balpha}_{\w}\right)^{\top}\bxi_\w^k \ge \h^\top_\w\y_\w^k , \forall \w\in\Omega, k \in K_{\w}^{(j)}   \: \label{Master2}\\
	\; \quad \;\, \W_\w \y_\w^k\; \geq \, \mathbf{b}_\w + \B_\w \bxi_\w^k - \T_\w \x,   \forall \w\in\Omega, k \in K_{\w}^{(j)}  \label{Master3}.
\end{IEEEeqnarray}

\noindent Where $K_{\w}^{(j)}=\{0,1,...,|\mathcal{E}_\w^{(j)}|\}$ is the set of indices that enumerates the set of scenarios $\mathcal{E}_\w^{(j)}=\{\bxi_{\omega}^{0}, \bxi_{\omega}^{1},...,\bxi_{\omega}^{|\mathcal{E}_\w^{(j)}|}\}$.

\noindent \textbf{3) Dantzig--Wolfe Procedure inner loop (updating $\mathcal{E}_\w^{(j)})$}: 
For each $\omega\in\Omega$, the counter $n$ is initialized with $n\leftarrow 1$, and the DWP is applied to the following constrained version of recourse problem \eqref{def_HMILP}--\eqref{def_HMILP2}:
\begin{IEEEeqnarray}{lrll}
	\underline{H}_{DR}                &(\x^{(j)},\w)=   \max_{p^k} &  \sum_{k\in K_{\w}^{(j)}} g(\x^{(j)},\bxi^k_{\omega},\w)\, p^k                \label{DWP1}   \\
	&\inmat{s.t.:}   &\sum_{k\in K_{\w}^{(j)}}  p^k = 1 &: \alpha_{0\w}^{(j)}  \label{DWP2} \\ 
	& \underline{\bmu}_{\w} \leq & \sum_{k\in K_{\w}^{(j)}} {\bxi^k_{\omega}} \,p^k \leq \overline{\bmu}_{\w} &:\underline{\balpha}_{\w}^{(j)},\overline{\balpha}_{\w}^{(j)},\;\;\;\;\,  \label{DWP3} 
\end{IEEEeqnarray}
where $\underline{H}_{DR}(\x^{(j)},\w)$ is a lower bound for ${H}_{DR}(\x^{(j)},\w)$. The DWP aims at improving the approximation for ${H}_{DR}(\x^{(j)},\w)$ by determining new columns for \eqref{DWP1}--\eqref{DWP3}. Thus, we construct an oracle subproblem to identify the scenario $\bxi_\w^*\in\Xi_{\omega}$ with the highest reduced cost, $\overline{c}_\w^*$: 
\begin{IEEEeqnarray}{rl'l}    
	(\overline{c}_\w^*,\bxi_\w^*) \leftarrow \max_{\bxi\in \Xi_{\omega}} & \Big\{ g(\x^{(j)},\bxi,\w) - \alpha_{0\w}^{(j)}  -  \left(\overline{\balpha}_{\w}^{(j)} - \underline{\balpha}_{\w}^{(j)}\right)^{\top}\bxi \Big\}.  \quad  \label{innerP}  
\end{IEEEeqnarray}

 To solve problem \eqref{innerP}, we replace function $g(\x^{(j)},\bxi,\w)$ with { the objective function of the dual problem associated with \eqref{compac}; that is $\left(\mathbf{b}_\w+\B_\w \bxi - \T_\w \x^{(j)} \right)^\top\bpi$, where $\bpi$ represents the dual vector for the constraints in \eqref{compac}. This modification to \eqref{innerP}} results in the following nonlinear problem:
\begin{IEEEeqnarray}{lll}  \label{nonlinearoracle}
    \max_{\substack{\bxi\in \Xi_{\omega}\\\bpi\geq 0}} & \,  \left(\mathbf{b}_\w+\B_\w \bxi - \T_\w \x^{(j)} \right)^\top\bpi  - \alpha_{0\w}^{(j)}  -  \left(\overline{\balpha}_{\w}^{(j)} - \underline{\balpha}_{\w}^{(j)}\right)^{\top}\bxi \,  \;\;\quad\label{dualobj}  \\
	 \;\inmat{s.t.:} &   \quad \W^{\top}_\w\bpi   \leq \h_\w. 
\end{IEEEeqnarray}
%

 Because $\bxi_\w^*$ belongs to the set of extreme points ($\mathcal{E}_\w$) of the box-like support set $\Xi_{\w}$\footnote{{The} optimal scenarios identified by problem \eqref{innerP} are extreme points of the polyhedral support set $\Xi_{\omega}$, as this problem is a maximization of a convex function within a polyhedral set \cite{bertsimas1997introduction}. { This fact aligns with the result derived in Section \ref{sec.finite}.}}, the value of each component $\xi[i]$ of the decision vector $\bxi$, in the optimal solution, is given either the upper ($\overline{\xi}[i]$) or lower ($\underline{\xi}[i]$) limit of the $i$-th dimension of $\Xi_{\w}$. Hence, an auxiliary binary variable $u[i]\in\{0,1\}$ is used to express each component in the optimal solution as follows: $\xi[i] = u[i]\overline{\xi}[i]+(1-u[i])\underline{\xi}[i]$. Under this transformation, the bilinear product $\bxi^\top\bpi$ is linearized by standard integer modeling techniques\footnote{{The linearization procedure requires the consideration of bounds for the components of the dual vector $\bpi$ associated with constraint \eqref{second_balance}. These bounds are set based on an iterative process that scales its initial value whenever one of the bounds are met. The initial value is set to the maximum imbalance cost. In our numerical results, we observed no particular numerical problem.}}, as done in \cite{moreira2016reliable,moreira2018}, and \cite{wei2016}.


\vspace{0.2 cm}
\noindent \textbf{4) Convergence of Dantzig--Wolfe inner loop:} At each iteration of the DWP inner loop, the short-term scenario $\bxi^*_\w$ determined by the oracle subproblem updates the set of extreme points, i.e., $\mathcal{E}_\w^{(j)} \leftarrow \mathcal{E}_\w^{(j)} \cup \{\bxi^*_\w$\}. 

{If $n$ equals $L$ or the DWP tolerance level is achieved, i.e., $\overline{c}_\w^* \le \epsilon$, the algorithm proceeds \textbf{to step 5}. Otherwise, it returns \textbf{to step 3} and $n\leftarrow n+1$}.

It is relevant to highlight that the  DWP is performed separately for each $\w \in\Omega$. Thus, $\forall\w \in\Omega$, the DWP identifies up to $L$ scenarios, improving from below, at each inner iteration, the approximation $\underline{H}_{DR}(\x^{(j)},\w)$ for ${H}_{DR}(\x^{(j)},\w)$.  


\vspace{0.2 cm}
\noindent \textbf{5) Convergence test for the main loop (assessing $UB^{(j)}$):} {The upper bound for the ECCG algorithm at iteration $j$ is given by $UB^{(j)}=\mathbf{c}^{\top}_{inv}\x^{(j)} + \sum_{\w\in\Omega} \rho_{\w} \overline{H}_{DR}(\x^{(j)},\w)$, where $\overline{H}_{DR}(\x^{(j)},\w) = \underline{H}_{DR}(\x^{(j)},\w) + \overline{c}_\w^*$. This result is demonstrated in \textbf{Theorem 1}, presented in Section \ref{teo1}.

If 
$1-LB^{(j)}/UB^{(j)} \le \varepsilon$, the algorithm terminates and $\x^{(j)}$ is the $\varepsilon$-optimal solution. Otherwise, for each $\omega\in\Omega$, the scenarios added by the DWP are ranked according to their contribution to $\underline{H}_{DR}$, i.e., the product $g(\x^{(j)},\bxi^k_{\omega},\w) {p^k}$. The $M \le L$ scenarios with highest ranks are added to the master problem for the next iteration. %
 This operation is represented in Fig. \ref{fig: ECCG} by $\mathcal{E}_\w^{(j)}\leftarrow best_M(\mathcal{E}_\w^{(j)}\backslash \mathcal{E}_\w^{(j-1)})$, where $\mathcal{E}_\w^{(-1)} = \emptyset$. Then, $\mathcal{E}_\w^{(j+1)}\leftarrow\mathcal{E}_\w^{(j)}$, update the iteration counter $j\leftarrow j+1$ and the algorithm returns \textbf{to step 2}. 
 
\textbf{\emph{Remark}:} Note that for a given pair $\{\x^{(j)}, \mathcal{E}_\w^{(j)}\}$, the proposed ECCG algorithm with  $M=L$ provides equal or greater lower bound for the iteration $j+1$ than that obtained by the CCG algorithm. This is because the CCG algorithm is a particular case that adds $M=L=1$ of the $M$ scenarios added by the general ECCG algorithm.

\subsection{Dantzig-Wolfe-like upper bound for the DRO{-TEP} model}
\label{teo1}

\vspace{0.2 cm}
\noindent \textbf{Theorem 1} (ECCG algorithm upper bound)
	\begin{enumerate}
	\item \textit{$UB^{(j)}=\mathbf{c}^{\top}_{inv}\x^{(j)} + \sum_{\w\in\Omega} \rho_{\w} \overline{H}_{DR}(\x^{(j)},\w)$ is an upper bound for the DRO problem \eqref{Pareto_DRO}, {i.e., $z^*_{DR} \le UB^{(j)}$}.}
	\vspace{0.1cm}
	\item \textit{$UB^{(j)}$ is $\epsilon$-tight for $\x^{(j)}$ when $L=+\infty$; that is, $UB^{(j)}\leq\mathbf{c}^{\top}_{inv}\x^{(j)} +H_{DR}(\x^{(j)},\w) + |\Omega|\epsilon$, particularly, for $\epsilon=0$, $UB^{(j)}=\mathbf{c}^{\top}_{inv}\x^{(j)} +H_{DR}(\x^{(j)},\w)$.}
	\end{enumerate}

%

  \vspace{0.2 cm}
The proof of \textbf{Theorem 1} is better presented using the result of \textbf{Proposition 1}: \textit{For each $\w\in\Omega$ and trial solution $\x^{(j)}$, $\overline{H}_{DR}(\x^{(j)},\w)=\underline{H}_{DR}(\x^{(j)},\w) + \overline{c}_\w^*$ is a local upper bound for the recourse function, i.e., $\overline{H}_{DR}(\x^{(j)},\w) \ge H_{DR}(\x^{(j)},\w)$ for all $\w\in\Omega$ and $\x^{(j)}$.} 

\noindent \textbf{Proof of Proposition 1}. Note that, for each $\w$, owing to optimality of problem \eqref{innerP}, 
\begin{IEEEeqnarray}{rl'l}    
	g(\x^{(j)},\bxi,\w) - \alpha_{0\w}^{(j)}  - \left(\overline{\balpha}_{\w}^{(j)} - \underline{\balpha}_{\w}^{(j)}\right)^{\top}\bxi \leq \overline{c}_\w^*, \quad \forall \bxi \in \Xi_{\omega}. \nonumber 
\end{IEEEeqnarray}
As the former inequality holds for each $\bxi$, it also holds on average for all probability measures in $\mathcal{D}_{\w}$; that is, $\forall P \in \mathcal{D}_{\w}$,
\begin{IEEEeqnarray}{rl'l}    
	\int_{\S_{\w}} g(\x^{(j)},\bxi,\w) dP - \alpha_{0\w}^{(j)} - \left(\overline{\balpha}_{\w}^{(j)} - \underline{\balpha}_{\w}^{(j)}\right)^{\top} \int_{\S_{\w}} \bxi dP \leq \overline{c}_\w^*. \nonumber 
\end{IEEEeqnarray}
Rearranging the terms and using the first-moment constraint in \eqref{eq. Def_D}, the following inequalities also hold for all $\forall P \in \mathcal{D}_{\w}$:
\begin{IEEEeqnarray}{rl'l}    
	\int_{\S_{\w}} g(\x^{(j)},\bxi,\w) dP\;\,&\leq \overline{c}_\w^* + \alpha_{0\w}^{(j)} + (\overline{\balpha}_{\w}^{(j)} - \underline{\balpha}_{\w}^{(j)})^{\top} \int_{\S_{\w}} \bxi dP \nonumber \\ 
	&\le \overline{c}_\w^* + \alpha_{0\w}^{(j)} + \overline{\balpha}_{\w}^{(j)\top} \overline{\bmu}_{\w} - \underline{\balpha}_{\w}^{(j)\top}\underline{\bmu}_{\w}. \nonumber 
\end{IEEEeqnarray}
Particularly, for the worst-case $P$ we have that $\int_{\S_{\w}} g(\x^{(j)},\bxi,\w) dP = H_{DR}(\x^{(j)},\w)$. Thus, we complete the proof of \textbf{Proposition 1} noting that, by strong duality (applied to problem \eqref{DWP1}--\eqref{DWP3}), the right-hand-side of the final inequality precisely meets our local upper bound. Therefore, 
\begin{IEEEeqnarray}{rl'l}    
	H_{DR}(\x^{(j)},\w)\;\,&\leq \overline{c}_\w^* + \alpha_{0\w}^{(j)} + \overline{\balpha}_{\w}^{(j)\top}\overline{\bmu}_{\w} - \underline{\balpha}_{\w}^{(j)\top}\underline{\bmu}_{\w} \nonumber\\
	&= \underline{H}_{DR}(\x^{(j)},\w) + \overline{c}_\w^*. \;\square \nonumber 
\end{IEEEeqnarray}


\noindent \textbf{Proof of Theorem 1}. The proof for item 1) is based on the result of \textbf{Proposition 1}. We know that
\begin{IEEEeqnarray}{rl'l}    
	\sum_{\w\in\Omega} \rho_{\w} H_{DR}(\x^{(j)},\w) \le \sum_{\w\in\Omega} \rho_\w \overline{H}_{DR}(\x^{(j)},\w), \nonumber 
\end{IEEEeqnarray}
as $\rho_\w\ge0$ for all $\w\in\Omega$. Hence, as $\x^{(j)}$ might not be the optimal solution, $\x^{*}$, 
\begin{IEEEeqnarray}{rl'l}    
	z^*_{DR}\,&=\mathbf{c}^{\top}_{inv}\x^{*} + \sum_{\w\in\Omega} \rho_{\w} H_{DR}(\x^{*},\w) \nonumber \\
	&\le \mathbf{c}^{\top}_{inv}\x^{(j)} + \sum_{\w\in\Omega} \rho_{\w} \overline{H}_{DR}(\x^{(j)},\w)=UB^{(j)}.  \nonumber 
\end{IEEEeqnarray}

As for the item 2), by \textbf{Proposition 1}, we have $\underline{H}_{DR}(\x^{(j)},\w)\leq H_{DR}(\x^{(j)},\w)\leq \overline{H}_{DR}(\x^{(j)},\w)$. Since the DWP converges for $L=+\infty$, then $\overline{c}_\w^*\leq\epsilon$ and  $\overline{H}_{DR}(\x^{(j)},\w) - H_{DR}(\x^{(j)},\w)\leq\epsilon$, $\forall\w\in\Omega$. Thus, $0\leq UB^{(j)}-\mathbf{c}^{\top}_{inv}\x^{(j)}-H_{DR}(\x^{(j)},\w)  \leq |\Omega|\epsilon$. $\square$


}

\section{Computational tests}
\label{sec.Results}

This section reports the results from computational experiments. {The first set of experiments is presented in Section \ref{subsec.Case}, where the solutions of DRO-TEP models are compared, in terms of costs and reliability, with the solutions of benchmark models, namely ARO-TEP and D-TEP}. {In the second set of experiments, in Section \ref{subsec.computational}, the computational capability of the proposed ECCG algorithm is compared with the} CCG algorithm and the FVA, for different uncertainty vector dimensions. 

 The computational experiments are based on a modified version of the IEEE 118-bus system, which encompasses 118 buses, 91 loads, 54 generators, 154 existing lines and 32 candidate lines for transmission expansion. Without loss of generality, we assume that the net demand ($\tilde{\bxi}$) is the only uncertain parameter. In the experiments, $\tilde{\bxi} = [\tilde{\xi_1},\tilde{\xi_2},\ldots,\tilde{\xi_d}]$ is modeled as the daily net load, composed of $d$ time periods, within a static TEP study. The uncertain net demand is decoupled into nodal demands by submatrix $\B_t$, as per \eqref{second_balance}. For expository purposes, in all experiments in this section, we used two distinct long-term scenarios, i.e., $\Omega = \{\omega_1, \omega_2\}$, by defining $\{[\underline{\bmu}_{\w}, \overline{\bmu}_{\w}]\}_{\w\in\Omega}$ and $\{\Xi_\w\}_{\w\in\Omega}$. In order to create diversity we stressed the system differently in each long-term scenario; that is, buses experiencing uncertain net demand in $\w_1$ and $\w_2$ are located in distinct regions of the network and are subjected to different parameters. In summary, $\w_1$ and $\w_2$ were constructed using two main principles: \begin{enumerate} \item $\w_1$ is much more volatile than $\w_2$; many more buses in $\w_1$ are subjected to uncertain net demand than in $\w_2$. \item $\w_2$ produces slightly higher costs than $\w_1$ for $\bxi\in[\underline{\bmu}_{\w}, \overline{\bmu}_{\w}]$ for the initial (not expanded) network; that is, $g(\mathbf{0},\bxi\in[\underline{\bmu}_{\w_2}, \overline{\bmu}_{\w_2}],\w_2)\geq g(\mathbf{0},\bxi\in[\underline{\bmu}_{\w_1}, \overline{\bmu}_{\w_1}],\w_1)$.\end{enumerate} For the sake of reproducibility, detailed input data used in the experiments of both Sections \ref{subsec.Case} and \ref{subsec.computational} can be downloaded from \cite{DROSystem}. All tests were run using Gurobi 7.0.2 under JuMP (Julia 0.5) on a Xeon E5-2680 processor at 2.5 GHz and 128 GB of RAM.

 %

\subsection{Assessment of the DRO-TEP model}
\label{subsec.Case}

We assessed instances of the proposed DRO-TEP model against instances of the benchmark models ARO-TEP and {D-TEP} in terms of cost and reliability. For quick reference, in this subsection, the instances of DRO-TEP, ARO-TEP and D-TEP models are denoted, respectively, by {DRO}($\rho_{\w_1}$,$\rho_{\w_2}$),  {ARO}($\rho_{\w_1}$,$\rho_{\w_2}$), and DET($\rho_{\w_1}$,$\rho_{\w_2}$), where $\rho_{\w_1}$ and $\rho_{\w_2}$ refer to the weights of corresponding long-term scenario. The specification of the instances of the benchmark models was performed by adapting the ambiguity set information as described in Section \ref{sec.particular}. In this experiment, the daily net demand was divided into six blocks of 4 hours ($d=6$) and {all instances were optimized to optimality.}

\begin{table*}[t]
	\centering \footnotesize
	\caption{TEP problems solutions and out-of-sample assessment}
	\begin{tabular}{p{1.35cm}|p{0.21cm}p{0.30cm}p{0.30cm}|p{0.33cm}p{0.33cm}p{0.33cm}p{0.33cm}|p{0.33cm}p{0.33cm}p{0.33cm}p{0.33cm}|p{0.33cm}p{0.33cm}p{0.33cm}p{0.33cm}|p{0.33cm}p{0.33cm}p{0.33cm}p{0.33cm}}
		\toprule
		& \multicolumn{3}{c|}{\multirow{3}{*}{TEP solutions}} & \multicolumn{16}{c}{\multirow{2}{*}{Out-of-sample simulations}} \\
		&&&&\multicolumn{1}{c}{}&\multicolumn{1}{c}{}&\multicolumn{1}{c}{}&\multicolumn{1}{c}{}&\multicolumn{1}{c}{}&\multicolumn{1}{c}{}&\multicolumn{1}{c}{}&\multicolumn{1}{c}{}&\multicolumn{1}{c}{}&\multicolumn{1}{c}{}&\multicolumn{1}{c}{}&\multicolumn{1}{c}{}&\multicolumn{1}{c}{}&\multicolumn{1}{c}{}&&\\
		& & &&\multicolumn{8}{c}{Low-variance distributions} & \multicolumn{8}{c}{High-variance distributions} \\
		\cline{2-20}\noalign{\vskip 0.85mm}
		&\multicolumn{1}{c}{\multirow{3}{*}{\shortstack{New\\lines}}}&  & & \multicolumn{4}{c|}{Normal distribution} & \multicolumn{4}{c|}{Beta distribution}  & \multicolumn{4}{c|}{Normal distribution} & \multicolumn{4}{c}{Beta distribution}\\
		&&\multicolumn{2}{c|}{Cost[$10^4$\$]}&\multicolumn{2}{c}{Cost[$10^4$\$]}&\multicolumn{2}{c|}{RI[\%]}&\multicolumn{2}{c}{Cost[$10^4$\$]}&\multicolumn{2}{c|}{RI[\%]}&\multicolumn{2}{c}{Cost[$10^4$\$]}&\multicolumn{2}{c|}{RI[\%]}&\multicolumn{2}{c}{Cost[$10^4$\$]}&\multicolumn{2}{c}{RI[\%]}\\
		\multicolumn{1}{c|}{\multirow{-6}{*}{Instance}}&&\multicolumn{1}{c}{Inv.}&\multicolumn{1}{c|}{Total}& $\w_1$ &$\w_2$& $\w_1$ &$\w_2$& $\w_1$ &$\w_2$& $\w_1$ &$\w_2$& $\;\w_1$ &$\w_2$& $\w_1$ &$\w_2$& $\w_1$ &$\w_2$& $\w_1$ &$\w_2$\\
		\cline{1-20}\noalign{\vskip 0.95mm}
		DRO(100,0)  & \multicolumn{1}{c}{9} & \multicolumn{1}{r}{19} & \multicolumn{1}{r|}{90}&61&67&\multicolumn{1}{r}{4}&\multicolumn{1}{r|}{1}& 55&61&\multicolumn{1}{r}{0}&\multicolumn{1}{r|}{0}&\multicolumn{1}{r}{69}&71&\multicolumn{1}{r}{8}&\multicolumn{1}{r|}{4} &57&66&\multicolumn{1}{r}{1} &\multicolumn{1}{r}{0} \\
		DRO(0,100)  & \multicolumn{1}{c}{3} & \multicolumn{1}{r}{7} & \multicolumn{1}{r|}{49}&85&48&\multicolumn{1}{r}{33}&\multicolumn{1}{r|}{1}&59&47&\multicolumn{1}{r}{11}&\multicolumn{1}{r|}{0}&\multicolumn{1}{r}{103}&50&\multicolumn{1}{r}{43}&\multicolumn{1}{r|}{2} &81&47&\multicolumn{1}{r}{34}&\multicolumn{1}{r}{0} \\
		DRO(50,50)  & \multicolumn{1}{c}{8} & \multicolumn{1}{r}{14} & \multicolumn{1}{r|}{76}&57&55&\multicolumn{1}{r}{5}&\multicolumn{1}{r|}{1} &50&54&\multicolumn{1}{r}{0}&\multicolumn{1}{r|}{0}&\multicolumn{1}{r}{65}&57&\multicolumn{1}{r}{9} &\multicolumn{1}{r|}{2} &52&54&\multicolumn{1}{r}{1}&\multicolumn{1}{r}{0}  \\
		ARO(50,50) & \multicolumn{1}{c}{10} & \multicolumn{1}{r}{22}  &\multicolumn{1}{r|}{100} &64&62&\multicolumn{1}{r}{4}&\multicolumn{1}{r|}{1} &58&61&\multicolumn{1}{r}{0}&\multicolumn{1}{r|}{0}&\multicolumn{1}{r}{71}&64&\multicolumn{1}{r}{8}&\multicolumn{1}{r|}{2}  &59&62&\multicolumn{1}{r}{1}&\multicolumn{1}{r}{0}\\
		DET(50,50) & \multicolumn{1}{c}{-} & \multicolumn{1}{r}{-} & \multicolumn{1}{r|}{38}    &81&64&\multicolumn{1}{r}{34}&\multicolumn{1}{r|}{23} &52&45&\multicolumn{1}{r}{12}&\multicolumn{1}{r|}{1}&\multicolumn{1}{r}{102}&76&\multicolumn{1}{r}{44}&\multicolumn{1}{r|}{34} &76&62&\multicolumn{1}{r}{35}&\multicolumn{1}{r}{21}\\
		\cline{1-20}\noalign{\vskip 0.95mm}
	\end{tabular}%
	\label{tab:TEPsolutions}%
\end{table*}%

{The experiment was conducted as follows. First, the TEP solution was computed for each instance. Then, the solutions were evaluated in out-of-sample simulations under variants of the Normal and Beta distributions for describing the net demand $\tilde{\bxi}$. Each out-of-sample simulation accounted for 10,000 samples of the short-term net demand, representing 5,000 simulated days for each $\w \in \Omega$. In these simulations, for each individual scenario realization, we used \eqref{second_fobj}--\eqref{eq.lastConst} to compute the operational dispatch. In total, 4 different distributions were used for each $\w \in \Omega$, resulting in 8 out-of-sample cases. As low-variance marginal distributions for $\tilde{\xi_t}$, we used: 1) a Normal distribution such that the probability of the event $[\tilde{\xi_t}\in \Xi_{tw} |\,\S_{\w}]$ is equal to 95\%, where $\Xi_{tw}$ is the projection of $\Xi_{w}$ for each period $t$; and 2) a symmetric Beta  distribution defined in $\Xi_{tw}$, with parameter $4.5$. As high-variance marginal distributions for $\tilde{\xi_t}$ we used: 1) a Normal distribution such that the probability of the event $[\tilde{\xi_t}\in \Xi_{tw} |\,\S_{\w}]$ is equal to 90\%; and 2) a symmetric Beta distribution defined in $\Xi_{tw}$, with parameter $1.5$. {Even though the net loads typically exhibit temporal correlation, for the sake of simplicity, we assumed independence across different periods.}

	The main outcomes of the case study are summarized in Table \ref{tab:TEPsolutions}, where column 1 identifies the instances. Columns 2--4 present TEP solutions, namely, number of invested lines (column 2), investment costs (column 3), and in-sample total costs (column 4) -- including the investment, operational, and imbalance costs. Therefore, columns 2 to 4 are classified as in-sample results. Columns 5--20 present out-of-sample results for selected distributions for $\tilde{\bxi}$ for both $\w_1$ and $\w_2$. The out-of-sample results are divided in 4 blocks, each block corresponding to one different distribution. For each block we present, for both $\w_1$ and $\w_2$, the total expected cost and the reliability index (RI). We define RI$_\w$ as the conditional probability (conditioned to long-term scenario $\w$) of experiencing a load shedding greater than 0.5\% of the system nominal demand within a given operative day. Despite the fact that the total cost already incorporates a term related to the dispatch infeasibility (imbalance costs), the RI provides valuable statistical (frequency) information regarding the operative reliability of the solutions.
	
	Regarding in-sample outcomes, the solution for DRO(100,0) is more conservative than that for DRO(0,100), determining more invested lines (9 vs. 3) at a higher investment cost (\$19x10$^4$ vs. \$7x10$^4$). This result is coherent since $\w_1$ was designed to be much more volatile than $\w_2$, which is confirmed by respective total in-sample costs presented in column 4. The solution for DRO(50,50), considering both $\w_1$ and $\w_2$, lies between DRO(100,0) and DRO(0,100) in terms of invested lines, investment costs and total cost. As expected, ARO(50,50) is the instance that presents the most conservative solution, encompassing 10 new lines at the cost of \$22x10$^4$. The solution for the instance DET(50,50), which solely relies on long-term averages (disregarding short-term variability), determined no investments and is the least expensive (\$38x10$^4$) in terms of total in-sample cost. Nevertheless, in-sample total costs are not comparable across different instances, since, by construction, they assume distinct distributions for $\tilde{\bxi}$.

As for the out-of-sample assessment, note that DRO(100,0) performs better than DRO(0,100) in terms of cost and RI for $\w_1$ across all distributions, whereas the opposite holds for $\w_2$. This result is coherent since the invested lines for these two instances were determined considering solely the corresponding long-term scenario, disregarding the outcomes of the other long-term scenario. In summary, when assessed under $\w_2$, DRO(100,0) presents reasonable RI, but high relative costs compared to other DRO-based instances. When assessed under $\w_1$, especially for high-variance distributions, DRO(0,100) presents very high RI and high costs. In contrast, DRO(50,50) performs well under all distributions in terms of costs and reliability. It is particularly interesting to highlight that DRO(50,50) presents 1) similar values for RI as those presented by DRO(0,100) under $\w_2$, and 2) lower costs than DRO(100,0) for all distributions, even under $\w_1$. The latter is explained by the fact that, as compared to DRO(50,50), the expensive investments of DRO(100,0), \$19x10$^4$, which were driven by the worst-case probability for the aggressive $\w_1$ only, do not pay off in terms of costs in the proposed out-of-sample analysis under more reasonable (unimodal) distributions.  
 
Regarding the benchmark instances, ARO(50,50) consistently presents low RI values under all distributions for both $\w_1$ and $\w_2$, as expected. Nevertheless, it is dominated in all cases, in terms of expected cost, by DRO(50,50), which also achieved almost the same values for RI, even under high-variance distributions. As for DET(50,50), it performs well under the low-variance Beta distribution for long-scenario $\w_2$. However, for all other cases, DET(50,50) is dominated by DRO(50,50) in terms of cost and reliability. It is relevant to highlight that the optimistic solution of DET(50,50) produces high costs as a consequence of the unacceptable risk levels observed in most cases, which rules out this approach under this experiment parameterization.

As the planner needs to decide the investment strategy under uncertainty of the distribution and long-term scenario, the authors advocate that, within  the  limitations  of  this  case  study, the DRO(50,50) instance produces the best tradeoff in terms of cost and RI among the tested instances. The authors acknowledge however that the extent of the benefits of the proposed approach can be influenced by the system topology, the premises for long-term scenarios, and  costs (investments, imbalance, etc.). Notwithstanding, the proposed DRO-TEP approach remains a relevant tool for planners in determining effective investment decisions featuring a consistent balance between expected cost and reliability in a setting where short-term operational uncertainty is ambiguous. 


\subsection{Assessment of the ECCG algorithm}
\label{subsec.computational}

{In this subsection, the performance of the proposed ECCG algorithm is assessed against those of the benchmark methods mentioned in Section \ref{sec.Methodology}, namely, the FVA applied in \cite{PozoPscc2018}, and the CCG algorithm applied in \cite{wei2016,zhao2017distributionally}}. All 3 methods (FVA, CCG algorithm, and ECCG algorithm) were applied to solve 4 different DRO(50,50) instances comprising {4, 5, 8, and 12 time periods ($d=4$, $5$, $8$ or $12$)}. 

In order to decide which results to report, we performed a sensitivity analysis on the parameters. As detailed in {S}ection \ref{sec.Methodology}, the specification of the ECCG algorithm involves 4 parameters / tolerances, namely, the maximum number of inner iterations $L$, the maximum number of scenarios included to the master problem per main iteration $M$, the DWP tolerance $\epsilon$, and global tolerance for the problem (main loop) $\varepsilon$. The global tolerance, which was set to $\varepsilon=1\%$, is the only parameter also applicable to FVA and CCG algorithm. We recall that the CCG algorithm is a particular case of the ECCG algorithm where $M=L=1$ and the tolerance $\epsilon$ is not applicable. As for the 3 other parameters of the ECCG method, our experiments indicated that when $L$ is higher than a certain threshold and $\epsilon$ is low, very good approximations for  $H_{DR}(\x,\w)$ are obtained. These approximations can be measured, for example, by the ratio $\overline{c}_\w^*/\underline{H}_{DR}(\x,\w)$ (the lower, the better). For the instance sizes ($d=4$, $5$, $8$ or $12$) that we investigated, increasing $L$ after a certain value almost did not impact the $M$ scenarios selected by the method. Thus, for the sake of brevity and due to space limitations, we decided to report results for fixed $L=20$ and $\epsilon=0.10\$$. In this setting, the sensitive parameter for the ECCG algorithm is $M$, which is reported for values ranging from 1 to 5. For quick reference, the instances of ECCG {are denoted by ECCG($M$)}. In this experiment, the Gurobi optimality gap was set to 0.5\% for the master problem of both decomposition methods (and to $\varepsilon=1\%$ for the FVA). Table \ref{tab:CPU_TIME} displays the computing times and number of iterations

It can be observed in Table \ref{tab:CPU_TIME} that the ECCG algorithm dominates the others, particularly, for higher values of $d$. Superior results were achieved for higher values of $M$ that required fewer iterations of the master problem. The benchmark methods are only practical for low-dimensional uncertainty vectors. The FVA, for instance, failed to converge in 18 hours for $d\geq8$. The CCG algorithm is also limited to low values of $d$ (failing to converge for $d=12$), mainly because it produces upper bounds that are not tight, thereby requiring additional iterations to approximate $H_{DR}(\x,\w)$. 

\begin{table}[t!]
	\centering\scriptsize
	\caption{Comparative CPU times (s) and number of iterations}
	\begin{tabular}{lrr|rr|rr|rr}
		\toprule
		\multirow{3}{0.40cm}{Method} & \multicolumn{8}{c}{Uncertainty Vector Dimension}\\
		\cline{2-9} \noalign{\vskip 0.95mm}
		& \multicolumn{2}{c|}{4} & \multicolumn{2}{c|}{5} &  \multicolumn{2}{c|}{8} &\multicolumn{2}{c}{12} \\
			&  \multicolumn{1}{c}{{t. }}& {Iter. } &   \multicolumn{1}{c}{{t. }}& {Iter. }  &   \multicolumn{1}{c}{{t. }}& {Iter. }  &   \multicolumn{1}{c}{{t. }}& {Iter. }  \\
		& \multicolumn{1}{c}{(s)} &  &  \multicolumn{1}{c}{(s)}  &  & \multicolumn{1}{c}{(s)}  &  & \multicolumn{1}{c}{(s)}  &  \\
		\midrule
		FVA	     &170&	-		&1,825	& -		&    \multicolumn{1}{c}{T}	& -	 &     \multicolumn{1}{c}{T}		& 	-	\\
		CCG 	 &289& 9 		&625	& 10    &8,171	& 15 &     \multicolumn{1}{c}{T}		&	 \multicolumn{1}{r}{T}	\\
		ECCG(1)   &111& 5  		&195   	& 5 	&1,262 	& 6	 &14,510	&11 	\\
		ECCG(2)   & 75& 3  		&119   	& 3 	&955 	& 4  &12,247	& 7 	\\
		ECCG(3)   & 87& 3   		&158   	& 3 	&811 	& 3  &8,503		& 5	 	\\
		ECCG(4)   & 56& 2 		&90     & 2 	&1,168 	& 3  &6,776		& 4	 	\\
		ECCG(5)   & 62& 2  		&98  	& 2 	&549 	& 2  &4,083 	& 3 	\\
		\bottomrule\\
		\multicolumn{9}{l}{T - Time limit of 18 hours exceeded without convergence.}
	\end{tabular}%
	\label{tab:CPU_TIME}%
\end{table}

{Table \ref{tab:CPU_TIME2} extends the results presented in Table \ref{tab:CPU_TIME} for $d=8$, detailing each iteration (reported in column 2) of both CCG and ECCG methods. Columns 3-5 provide the computing times for the master problem, the DWP inner loop, and the whole iteration. {Columns 5-7 present}, respectively, the processing times of a given iteration, the cumulative time, and the gap of the ECCG algorithm. 
	{Columns 8-9 display} the number of inner loop iterations for $\w_1$ and $\w_2$.}

\begin{table}[t!]
	
	\centering \scriptsize
	\caption{{Detailed iteration data for $d$=8 }}
	\begin{tabular}{c|r|rrrr|r|rr}
		\toprule
		\rowcolor[HTML]{FFFFFF} 
		\cellcolor[HTML]{FFFFFF}  & \cellcolor[HTML]{FFFFFF} & \multicolumn{4}{c|}{\cellcolor[HTML]{FFFFFF}Time} & \multicolumn{1}{c|}{\cellcolor[HTML]{FFFFFF} GAP} & \multicolumn{2}{r}{\cellcolor[HTML]{FFFFFF}Inner iter.} \\
		\rowcolor[HTML]{FFFFFF} 
		\cellcolor[HTML]{FFFFFF}  & \cellcolor[HTML]{FFFFFF} & \multicolumn{4}{c|}{\cellcolor[HTML]{FFFFFF}(s)} & \multicolumn{1}{c|}{\cellcolor[HTML]{FFFFFF} (\%)} & \multicolumn{2}{c}{\cellcolor[HTML]{FFFFFF}(\#)} \\
		\rowcolor[HTML]{FFFFFF} 
		\multicolumn{1}{c|}{\multirow{-3}{*}{\cellcolor[HTML]{FFFFFF}Method}} & 	\multicolumn{1}{c|}{\multirow{-3}{*}{Iter.}} & \multicolumn{1}{c}{Master} & \multicolumn{1}{c}{\cellcolor[HTML]{FFFFFF}DWP} & \multicolumn{1}{l}{\cellcolor[HTML]{FFFFFF}Iter.} & \multicolumn{1}{r|}{\cellcolor[HTML]{FFFFFF}Accum.} & \multicolumn{1}{r|}{\cellcolor[HTML]{FFFFFF}} & \multicolumn{1}{r}{\cellcolor[HTML]{FFFFFF}$\w_1$} & \multicolumn{1}{r}{\cellcolor[HTML]{FFFFFF}$\w_2$} \\
		\cline{1-9}\noalign{\vskip 0.80mm}
		\rowcolor[HTML]{EFEFEF}&1& 13 & 7 & 20 & 20 & 85.1 & 1 & 1 \\
		\rowcolor[HTML]{EFEFEF}&2& 29 & 9 & 38 & 58 & 57.2 & 1 & 1 \\
		\rowcolor[HTML]{EFEFEF}&3& 47 & 9 & 56 & 114 & 56.1 & 1 & 1 \\
		\rowcolor[HTML]{EFEFEF}&4& 99 & 9 & 108 & 222 & 55.0 & 1 & 1 \\
		\rowcolor[HTML]{EFEFEF}&5& 142 & 11 & 153 & 375 & 45.8 & 1 & 1 \\
		\rowcolor[HTML]{EFEFEF}&6& 148 & 12 & 160 & 535 & 45.3 & 1 & 1 \\
		\rowcolor[HTML]{EFEFEF}& 7 & 244 & 11 & 255 & 790 & 44.8 & 1 & 1 \\
		\rowcolor[HTML]{EFEFEF}& 8 & 365 & 13 & 378 & 1,168 & 40.0 & 1 & 1 \\
		\rowcolor[HTML]{EFEFEF}  & 9 & 343 & 12 & 355 & 1,523 & 40.0 & 1 & 1 \\
		\rowcolor[HTML]{EFEFEF}  & 10 & 473 & 11 & 484 & 2,007 & 39.8 & 1 & 1 \\
		\rowcolor[HTML]{EFEFEF}  & 11 & 760 & 13 & 773 & 2,780 & 35.4 & 1 & 1 \\
		\rowcolor[HTML]{EFEFEF}  & 12 & 711 & 14 & 725 & 3,505 & 10.0 & 1 & 1 \\
		\rowcolor[HTML]{EFEFEF}  & 13 & 1,319 & 18 & 1,337 & 4,842 & 7.4 & 1 & 1 \\
		\rowcolor[HTML]{EFEFEF}  & 14 & 1,671 & 17 & 1,688 & 6,530 & 1.4 & 1 & 1 \\
		\rowcolor[HTML]{EFEFEF} 
		{\multirow{-15}{*}{CCG}} & 15 & 1,623 & 18 & 1,641 & 8,171 & 0.7 & 1 & 1 \\
		\cline{1-9}\noalign{\vskip 0.80mm}
		\rowcolor[HTML]{FFFFFF} 	& 1 & 13 & 81 & 94 & 94 & 60.8 & 16 & 20 \\
		\rowcolor[HTML]{FFFFFF} 	& 2 & 23 & 124 & 147 & 241 & 47.8 & 19 & 20 \\
		\rowcolor[HTML]{FFFFFF} 	& 3 & 60 & 128 & 188 & 429 & 41.3 & 19 & 20 \\
		\rowcolor[HTML]{FFFFFF} 	& 4 & 110 & 134 & 244 & 673 & 1.6 & 14 & 18 \\
		\rowcolor[HTML]{FFFFFF} 	& 5 & 147 & 119 & 266 & 939 & 1.0 & 13 & 14 \\
		\rowcolor[HTML]{FFFFFF} 	\multirow{-6}{*}{ECCG(1)} & 6 & 226 & 97 & 323 & 1,262 & 0.8 & 12 & 8 \\
		\cline{1-9}\noalign{\vskip 0.80mm}
		\rowcolor[HTML]{EFEFEF} 
		\cellcolor[HTML]{EFEFEF} & 1 & 13 & 82 & 95 & 95 & 60.8 & 16 & 20 \\
		\rowcolor[HTML]{EFEFEF} 
		\cellcolor[HTML]{EFEFEF} & 2 & 48 & 142 & 190 & 285 & 33.6 & 20 & 20 \\
		\rowcolor[HTML]{EFEFEF} 
		\cellcolor[HTML]{EFEFEF} & 3 & 169 & 134 & 303 & 588 & 10.4 & 18 & 17 \\
		\rowcolor[HTML]{EFEFEF} 
		\multirow{-4}{*}{\cellcolor[HTML]{EFEFEF}ECCG(2)} & 4 & 289 & 78 & 367 & 955 & 0.6 & 5 & 12 \\
		\cline{1-9}\noalign{\vskip 0.80mm}
		\rowcolor[HTML]{FFFFFF} 	& 1 & 14 & 83 & 97 & 97 & 60.8 & 16 & 20 \\
		\rowcolor[HTML]{FFFFFF} 	& 2 & 101 & 118 & 219 & 316 & 16.2 & 14 & 20 \\
		\rowcolor[HTML]{FFFFFF} 	\multirow{-3}{*}{ECCG(3)} & 3 & 386 & 109 & 495 & 811 & 0.7 & 10 & 14 \\
		\cline{1-9}\noalign{\vskip 0.80mm}
		\rowcolor[HTML]{EFEFEF} 
		\cellcolor[HTML]{EFEFEF} & 1 & 14 & 83 & 97 & 97 & 60.8 & 16 & 20 \\
		\rowcolor[HTML]{EFEFEF} 
		\cellcolor[HTML]{EFEFEF} & 2 & 151 & 153 & 304 & 401 & 3.9 & 14 & 20 \\
		\rowcolor[HTML]{EFEFEF} 
		\multirow{-3}{*}{\cellcolor[HTML]{EFEFEF}ECCG(4)} & 3 & 685 & 82 & 767 & 1,168 & 0.2 & 9 & 7 \\
		\cline{1-9}\noalign{\vskip 0.80mm}
		\rowcolor[HTML]{FFFFFF} 
		\cellcolor[HTML]{FFFFFF} & 1 & 14 & 82 & 96 & 96 & 60.8 & 16 & 20 \\
		\rowcolor[HTML]{FFFFFF} 
		\multirow{-2}{*}{\cellcolor[HTML]{FFFFFF}ECCG(5)} & 2 & 321 & 132 & 453 & 549 & 0.7 & 14 & 17\\
		\cline{1-9}\noalign{\vskip 0.95mm}
	\end{tabular}
	\label{tab:CPU_TIME2}%
\end{table}

 {Table \ref{tab:CPU_TIME2} results show that the time to solve the master problem grows across iterations with a much higher rate than the time to solve the subproblems for all cases. Therefore, it is clear that the master problem is the bottleneck for both decomposition methods. The key advantage of the ECCG algorithm over the CCG algorithm is that it requires less iterations of the master problem because of two main factors: 1) tighter upper bounds due to the DWP, and 2) tighter lower bounds due to the selection of up to $M$ scenarios with high contribution to the recourse function to be added to the master problem at each iteration. The effect of point 1 is observed by comparing the GAP after the first iteration of the CCG, 85.1\%, with the GAP for any ECCG method, 60.8\%. As the first stage solution is the same for all methods, the difference is solely due to the better approximation of the recourse function. The effect of point 2 is detected by the improved performance with higher values of $M$. 
 	



\section{Conclusions}
\label{sec: conclusions}

{In this work we {have} proposed a multi-scale distributionally robust transmission expansion {planning} model. {We have} introduced the concept of multiple conditional ambiguity sets to account for the information of long-term studies conducted by experts in current industry practices. Due to intractability issues associated with DRO{-}TEP models, an enhanced{-}column-and-constraint-generation algorithm, providing better approximations of the {recourse} function and tighter bounds, was devised. Results for the IEEE 118-bus system have shown that, in comparison to existing methods, the proposed DRO-TEP model is effective in producing a consistent tradeoff between cost and reliability in out-of-sample tests. {As for the computational aspects,} the ECCG algorithm has significantly {outperformed} the CCG algorithm. Notwithstanding, the proposed methodology is not exhaustively developed in this work. The foundations of ECCG rely on finding new scenarios for the recourse function, then recomputing the dual variables iteratively. Extending this rationale for different ambiguity sets is an avenue for research. Likewise, devising metrics to select the best $M$ scenarios for the master problem and studying procedures to drop old cuts are interesting {research themes}}. 

\bibliographystyle{IEEEtran}
\bibliography{References}

\end{document}